\magnification=\magstep1
\baselineskip=15pt
\input psfig.sty
\centerline{\bf FOLIATIONS ON COMPLEX PROJECTIVE SURFACES}
\vskip 2truecm
\centerline{\bf Marco Brunella}
\vskip 4truecm
\par
In this text we shall review the classification of foliations on complex
projective surfaces according to their Kodaira dimension, following McQuillan's
seminal paper [MQ1] with some complements and variations given by [Br1] and
[Br2]. Most of the proofs will be only sketched, and the text should be
considered as guidelines to the above works (and related ones), with no
exhaustivity nor selfcontainedness pretention. There are no new results, but
some old results are presented in a new way, by adopting systematically an
orbifold point of view.
\vskip 1truecm
\par
{\bf Contents}
\par
1. Basic definitions
\par
2. Basic formulae
\par
3. Singularities
\par
4. Nef models
\par
5. Numerically trivial foliations
\par
6. Kodaira dimension
\par
7. Riccati and Turbulent foliations
\par
8. Poincar\'e metric
\par
9. Hilbert modular foliations
\par
10. K\"ahler surfaces

\vfill\eject\par
{\bf 1. Basic definitions}
\vglue 1truecm\par

Let $X$ be a smooth complex surface (the smoothness assumption will be soon
relaxed). A {\bf foliation} ${\cal F}$ on $X$ is given by an open covering $\{
U_j\}$ of $X$ and holomorphic vector fields $v_j\in H^0(U_j,\Theta_X)$ with 
isolated zeroes such that
$$v_i=g_{ij}v_j \qquad {\rm on} \qquad U_i\cap U_j$$
for some nonvanishing holomorphic functions $g_{ij}\in H^0(U_i\cap U_j, {\cal
O}_X^*)$. This allows to glue together the local orbits of the vector fields $\{
v_j\}$ to obtain the {\bf leaves} of ${\cal F}$. The {\bf singular set}
$Sing({\cal F})$ of ${\cal F}$ is the discrete subset of $X$ defined by
$Sing({\cal F})\cap U_j=\{$ zeroes of $v_j\}$. 
\par
The functions $\{ g_{ij}\}$ form a multiplicative cocycle and define a
holomorphic line bundle $K_{\cal F}$, called {\bf canonical bundle} of ${\cal
F}$. The relations $v_i=g_{ij}v_j$ allow to construct a global holomorphic
section $s\in H^0(X,K_{\cal F}\otimes\Theta_X)$, vanishing only on $Sing({\cal
F})$. If $T_{\cal F}$ denotes the dual of $K_{\cal F}$ (the {\bf tangent bundle}
of ${\cal F}$) we therefore have an exact sequence of sheaves induced by $s$
$$0\longrightarrow T_{\cal F}\longrightarrow \Theta_X \longrightarrow 
{\cal I}_Z\cdot N_{\cal F} \longrightarrow 0$$
for a suitable line bundle $N_{\cal F}$ (the {\bf normal bundle} of ${\cal F}$)
and a suitable ideal sheaf ${\cal I}_Z$ supported on $Sing({\cal F})$ (see e.g.
[Fri, Chapter 2]).
\par
In a dual way, the foliation ${\cal F}$ can be also defined by (kernels of)
holomorphic 1-forms with isolated zeroes $\omega_j\in H^0(U_j,\Omega_X^1)$
satisfying
$$\omega_i =f_{ij}\omega_j \qquad {\rm on} \qquad U_i\cap U_j$$
for a suitable cocycle $\{ f_{ij}\}$. This cocycle actually defines the normal
bundle $N_{\cal F}$, and if $N_{\cal F}^*$ denotes its dual (the {\bf conormal
bundle} of ${\cal F}$) we obtain an exact sequence of sheaves
$$0\longrightarrow N_{\cal F}^*\longrightarrow \Omega_X^1 \longrightarrow 
{\cal I}_Z\cdot K_{\cal F} \longrightarrow 0$$
which is the dual of the previous one.
\par
The canonical and the conormal bundle of ${\cal F}$ are related by the formula
$$K_X=K_{\cal F}\otimes N_{\cal F}^*$$
where $K_X$ ($=\Omega_X^2$) is the canonical bundle of the surface $X$. Indeed,
the choice of a local nonvanishing holomorphic 2-form induces, by contraction, a
local isomorphism between vector fields and 1-forms generating ${\cal F}$.
\par
Let us consider now the case of a singular surface $X$. In fact, we shall need
only the case where $X$ has only {\it cyclic quotient singularities}, i.e.
around each $p\in Sing(X)$ the surface is of the type ${\bf B}^2/\Gamma_{k,h}$
where ${\bf B}^2$ is the unit ball in ${\bf C}^2$ and $\Gamma_{k,h}$ is the
cyclic group of order $k$ generated by $(z,w)\mapsto (e^{{2\pi i\over k}}z,
e^{{2\pi i\over k}h}w)$, for suitable coprime positive integers $k,h$ with $0 <
h < k$. Equivalently [BPV, pages 80-85], such a singularity arises from the 
contraction of a Hirzebruch--Jung string on a smooth surface (a chain of 
rational curves, each one of selfintersection $\le -2$).
\par
Because cyclic quotient singularities are normal, a foliation ${\cal F}$ on such
a surface can be simply defined as a foliation on $X\setminus Sing(X)$. If $p\in
Sing(X)$ and $U\simeq {\bf B}^2/\Gamma$ is a neighbourhood of $p$, then ${\cal
F}\vert_{U\setminus\{ p\} }$ can be lifted to ${\bf B}^2\setminus\{ 0\}$, and on
that covering the foliation can be defined by a holomorphic vector field which
holomorphically extends to $0$. We shall say that ${\cal F}$ is {\bf
nonsingular} at $p$ if that extension is nonvanishing at $0$. In other words, if
${\cal F}$ is not singular at $p$ then on $U\simeq {\bf B}^2/\Gamma$ the
foliation is the quotient of the vertical or horizontal foliation on ${\bf B}^2$
(up to an equivariant biholomorphism). Even if many things can be done in
greater generality, we shall always assume that ${\cal F}$ is not singular at
singular points of $X$:
$$Sing({\cal F})\cap Sing(X)=\emptyset .$$
\par
In this context, leaves of ${\cal F}$ are defined as in the smooth case, by
glueing together local leaves through nonsingular points of ${\cal F}$. The only
difference is that if $p\in Sing(X)$ has order $k$ then the local leaf of ${\cal
F}$ through $p$ is an {\it orbifold} in which $p$ is affected by the
multiplicity $k$. Indeed, this local leaf on ${\bf B}^2/\Gamma$ is the quotient
of a disc ${\bf D}\subset{\bf B}^2$ by a $k$-cyclic group. Thus, leaves of
${\cal F}$ are {\it orbifolds} injectively immersed (in orbifold's sense) in
$X\setminus Sing({\cal F})$, and giving a partition of $X\setminus 
Sing({\cal F})$ into disjoint subsets.
\par
Given a foliation ${\cal F}$ on a surface $X$ with (at most) cyclic quotient
singularities, we can still define its canonical sheaf $K_{\cal F}$, for
instance by taking the direct image under the inclusion $X\setminus Sing(X) \to
X$. However, it is no more a genuine line bundle (an element of $Pic(X)$) but
only a ${\bf Q}$-bundle (an element of $Pic(X)\otimes{\bf Q}$). Indeed, if $p\in
Sing(X)$ has order $k$ then $K_{\cal F}$ is not locally free at $p$, but its
$k$-power $K_{\cal F}^{\otimes k}$ is, and moreover $k$ is the minimal positive
integer with that property (the vector field ${\partial\over\partial w}$ on
${\bf B}^2$ is not $\Gamma_{k,h}$-invariant, but its power 
$({\partial\over\partial w})^{\otimes k}$ is). Anyway, for many things 
${\bf Q}$-bundles are as good as line bundles. 
\par
Similar considerations hold also for $T_{\cal F}$, $N_{\cal F}$, $N_{\cal F}^*$,
and of course $K_X$. We still have the equality $K_X=K_{\cal F}\otimes 
N_{\cal F}^*$, as sheaves or ${\bf Q}$-bundles.

\vskip 2truecm\par
{\bf 2. Basic formulae}
\vglue 1truecm\par

We shall need many times some elementary formulae which compute the degree of
$K_{\cal F}$ over compact curves in $X$.
\par
As before, let $X$ be a surface with at most cyclic quotient singularities, and
let ${\cal F}$ be a foliation on $X$, with $Sing({\cal F})\cap
Sing(X)=\emptyset$. 
Let $C\subset X$ be a compact connected (possibly singular) curve, 
and suppose that each irreducible component of $C$ is not invariant by ${\cal
F}$. For every $p\in C$ we can define an index $tang({\cal F},C,p)$ which
measure the tangency order of ${\cal F}$ with $C$ at $p$ (and thus which is 0
for a generic $p\in C$, where we have transversality).
\par
If $p\not\in Sing(X)$ the definition is the following: we take a local equation
$f$ of $C$ at $p$, a local holomorphic vector field $v$ generating ${\cal F}$
around $p$, and we set
$$tang({\cal F},C,p) = \dim_{\bf C}{{\cal O}_p\over <f,v(f)>}$$
where ${\cal O}_p$ is the local algebra of $X$ at $p$ (germs of holomorphic
functions), $v(f)$ is the Lie derivative of $f$ along $v$, $<f,v(f)>$ is the
ideal in ${\cal O}_p$ generated by $f$ and $v(f)$. This index is finite, because
$C$ is not ${\cal F}$-invariant; it is a nonnegative integer; it is 0 iff
$p\not\in Sing({\cal F})$ and ${\cal F}$ is transverse to $C$ at $p$.
\par
If $p\in Sing(X)$, we take a neighbourhood $U\simeq{\bf B}^2/\Gamma$ and we lift
${\cal F}\vert_U$ and $C\cap U$ on ${\bf B}^2$. We compute the index on ${\bf
B}^2$ (at 0) and we divide the result by $k$, the order of the singularity. This
is, by definition, $tang({\cal F},C,p)$. It is a nonnegative rational number,
and it is 0 iff ${\cal F}$ is transverse to $C$ at $p$ (in the sense that the
lift of ${\cal F}$ on ${\bf B}^2$ is transverse to the lift of $C$ at 0).
\par
We also have at our disposal the selfintersection $C\cdot C$ of $C$ in $X$: it 
is the degree on $C$ of the ${\bf Q}$-bundle ${\cal O}_X(C)$, and it is a 
rational number, possibly a noninteger one if $C$ passes through $Sing(X)$ (see
e.g. [Sak] for the intersection theory on normal surfaces). 
We then have the following formula for the degree of the ${\bf Q}$-bundle 
$K_{\cal F}$ on $C$:
$$K_{\cal F}\cdot C = -C\cdot C + tang({\cal F},C)$$
where $tang({\cal F},C)=\sum_{p\in C}tang({\cal F},C,p)$ (a finite sum). This is
proved in [Br1, page 23] when $X$ is smooth, but the general case can be handled
along the same lines. The most important consequence of this formula is the
inequality
$$(K_{\cal F}+C)\cdot C \ge 0$$
where the equality is realized if and only if ${\cal F}$ is everywhere
transverse to $C$.
\par
Let us consider now the case where each irreducible component of $C$ is ${\cal
F}$-invariant. Again, we shall define an index $Z({\cal F},C,p)$ for every $p\in
C$. If $p\not\in Sing({\cal F})$ we simply set $Z({\cal F},C,p)=0$. If $p\in
Sing({\cal F})$ (thus, in particular, $p\not\in Sing(X)$) then we choose a local
equation $f$ of $C$ and a local 1-form $\omega$ generating ${\cal F}$. We can
factorize around $p$ (see e.g. [Suw, Chapter V])
$$g\omega = hdf + f\eta \qquad (i.e.\ {\omega\over f}={h\over g}{df\over f} +
{1\over g}\eta )$$
where $\eta$ is a holomorphic 1-form, $g$ and $h$ are holomorphic functions, and
$g$ and $h$ do not vanish identically on each local branch of $C$ at $p$. Thus
$h/g$ is a meromorphic function on each local branch of $C$ at $p$
(the residue of ${\omega\over f}$ along the branch), not
identically zero nor infinity. It has a well defined vanishing or polar order at
$p$ on each local branch of $C$ at $p$, and the sum of these vanishing or polar
orders is, by definition, the index $Z({\cal F},C,p)$. 
\par
We refer to [Br1, page 24] for a more detailed discussion of this index. Here we
just recall that: i) the index is an integer number, but it may be negative,
when $C$ is a so-called ``dicritical separatrix at $p$''; ii) if $C$ is smooth
at $p$, the index equals the vanishing order of $v\vert_C$ at $p$, $v$ being a
local holomorphic vector field generating ${\cal F}$, and so it is a positive
number; iii) if $C$ has a normal crossing singularity at $p$, with local
branches $C_1$ and $C_2$, then $Z({\cal F},C,p)=Z({\cal F},C_1,p) + 
Z({\cal F},C_2,p)-2$, which is nonnegative.
\par
As in [Br1, page 25] we then have the formula:
$$K_{\cal F}\cdot C = -\chi_{orb}(C) + Z({\cal F},C)$$
where $Z({\cal F},C)=\sum_{p\in C} Z({\cal F},C,p)$. Here $\chi_{orb}(C)$
denotes the orbifold-arithmetic Euler characteristic of $C$, which may be
defined by the adjunction formula
$$\chi_{orb}(C)+C\cdot C = -K_X\cdot C$$
(each $p\in C\cap Sing(X)$ is affected by a multiplicity $k_p$, the order of the
singularity; if $C$ is smooth then $\chi_{orb}(C) = \chi_{top}(C) + \sum
{1-k_p\over k_p}$, the sum being over all $p\in C\cap Sing(X)$ and
$\chi_{top}(C)=2-2genus(C)$ being the topological Euler characteristic; if $C$
has singularities then one has the usual correction terms arising from smoothing
them [BPV, page 68]).
\par
Another very useful formula is Camacho--Sad formula 
[C-S] [Suw] [Br1, Chapter 3].
Let us consider again a compact connected curve $C$ invariant by ${\cal F}$. 
For each $p\in Sing({\cal F})\cap C$ we take a local factorization 
$g\omega =hdf+f\eta$
as above. On each local branch of $C$ at $p$ the meromorphic 1-form 
$-{1\over h}\eta$ has a well defined residue at $p$. The sum of these residues,
over all local branches of $C$ at $p$, is by definition the Camacho--Sad index
$CS({\cal F},C,p)$. Then we have [Br1, page 37]
$$C\cdot C = CS({\cal F},C) = \sum_{p\in Sing({\cal F})\cap C}
CS({\cal F},C,p).$$
We refer to [Br1, Chapter 3] or [Suw, Chapter V] for a more detailed discussion
of these Camacho--Sad index and formula.

\vskip 2truecm\par
{\bf 3. Singularities}
\vglue 1truecm\par
 
In order to simplify some statements and some proofs, and also to get a more
coherent theory, we will need some assumptions on the singularities of the
foliation [Br1, Chapter 1] [C-S] [M-M] [MR1] [MR2].
\par 
A singularity $p\in Sing({\cal F})$ is called {\bf reduced} (in Seidenberg's
sense) if ${\cal F}$ around $p$ is generated by a vector field $v$ whose linear
part $(Dv)_p$ has eigenvalues $\lambda_1 ,\lambda_2$ such that either $\lambda_1
\not= 0 \not= \lambda_2$ and $\lambda_1/\lambda_2\not\in {\bf Q}^+$, or
$\lambda_1\not= 0=\lambda_2$. In the former case $p$ is called {\bf
nondegenerate}, in the latter case a {\bf saddle-node}.
\par
A fundamental result of Seidenberg [M-M] [C-S] [Br1, page 13] says that, given 
any $p\in Sing({\cal F})$, we can perform a sequence of blow-ups based at $p$,
$\widetilde X\buildrel\pi\over\to X$, such that the lifted foliation 
$\widetilde{\cal F}$ on $\widetilde X$ has only reduced singularities on the 
exceptional divisor $\pi^{-1}(p)$.
Thus, from a bimeromorphic point of view, we can work without
loss of generality with {\bf reduced foliations}, i.e. foliations all of whose
singularities are reduced.
\par
Note also that the class of reduced singularities is stable by blow-ups: the
blow-up of a reduced singularity produces, on the exceptional divisor, only
reduced singularities. This remains true even if we blow up a regular point
$p\in X\setminus\{ Sing({\cal F})\cup Sing(X)\}$, or if we take the resolution
of a point $p\in Sing(X)$.

\centerline{\hbox{\psfig{figure=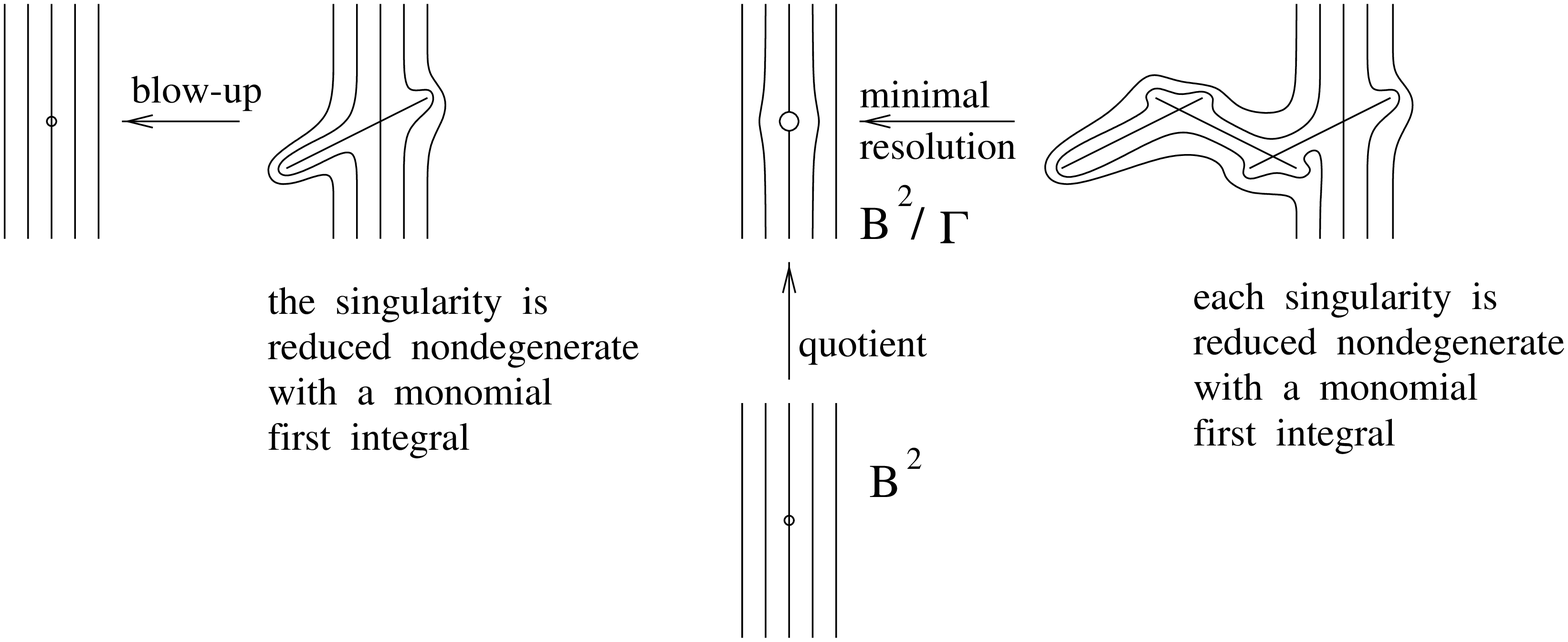,height=6truecm}}}

\par
Let us compute the indices $Z$ and $CS$ of the previous section, in the case of
reduced singularities. Suppose firstly that $p\in Sing({\cal F})$ is reduced and
nondegenerate. Then, in suitable local coordinates $(z,w)$ centered at $p$,
${\cal F}$ is generated by a vector field of the form [M-M] [Br1]
$$v=z{\partial\over\partial z} + \lambda w(1+...){\partial\over\partial w}$$
with $\lambda\not= 0$, $\lambda\not\in{\bf Q}^+$. The axis $\{ z=0\}$ and 
$\{ w=0\}$ are ${\cal F}$-invariant, and they are the only local curves through
$p$ invariant by ${\cal F}$. We easily find
$$Z({\cal F}, \{ z=0\} , 0) = Z({\cal F}, \{ w=0\} , 0) = 1 \qquad 
Z({\cal F}, \{ zw=0\} , 0) = 0$$
and
$$CS({\cal F}, \{ w=0\} , 0) = \lambda \qquad 
CS({\cal F}, \{ z=0\} , 0) = {1\over\lambda} \qquad 
CS({\cal F}, \{ zw=0\} , 0) = \lambda + {1\over\lambda} + 2 .$$
Suppose now that $p\in Sing({\cal F})$ is a saddle-node. In suitable local
coordinates $(z,w)$, ${\cal F}$ is generated by a vector field of the type
[M-M] [Br1]
$$v=[z(1+\nu w^k)+wF(z,w)]{\partial\over\partial z} +
w^{k+1}{\partial\over\partial w}$$
with $k\in{\bf N}^+$, $\nu\in{\bf C}$, and $F$ vanishes at $(0,0)$ up to order
$k$. The axis $\{ w=0\}$ (called {\bf strong separatrix}) is ${\cal
F}$-invariant, and we find
$$Z({\cal F},\{ w=0\} ,0)=1$$
$$CS({\cal F},\{ w=0\} ,0)=0 .$$
There exists {\it at most} one more local curve through $p$ and ${\cal
F}$-invariant. If it exists, it is smooth and transverse to the strong
separatrix, and it is called {\bf weak separatrix}. Its indices are given by
$Z=k+1$ and $CS=\nu$.
\par
Let us observe the following important consequence of the above discussion:
if ${\cal F}$ is a reduced foliation and $C$ is a ${\cal F}$-invariant curve,
then $C$ has {\it at most normal crossings} as singularities.
\par
We shall need also some information on the holonomy of the separatrices of
reduced singularities [Br1, pages 10-12]. If $p\in Sing({\cal F})$ is
nondegenerate, generated by a vector field $v$ as above, then the holonomy along
a loop in $\{ w=0\}$ has the form $h(w)=e^{2\pi i\lambda}w+...$. By a result of
Mattei and Moussu [M-M] [MR2], this holonomy characterizes the analytic 
conjugacy class of ${\cal F}$ around $p$. In particular, if $h$ is periodic of 
period $m$ (thus conjugate to $w\mapsto e^{{2\pi i n\over m}}w$ for some $n$) 
then ${\cal F}$ around $p$ is conjugate to the foliation generated by 
$mz{\partial\over\partial z} + lw{\partial\over\partial w}$, 
for some $l=n\ {\rm mod}\ m<0$ ($\lambda =l/m$).
Such a foliation has a local holomorphic first integral, given by $z^lw^m$.
\par
If $p\in Sing({\cal F})$ is a saddle-node, generated by a vector field $v$ as
above, then the holonomy along a loop in $\{ w=0\}$ has the form
$h(w)=w+w^{k+1}+...$, in particular it is never periodic. By a result of
Martinet and Ramis [MR1], this holonomy characterizes the analytic conjugacy 
class of ${\cal F}$ around $p$.
\par
Let us finally remark that if $p\in Sing(X)$ and $\gamma$ is a loop around $p$
in the leaf $L_p$ of ${\cal F}$ through $p$, then the holonomy of ${\cal F}$
along $\gamma$ is {\it not} the identity, but it is a periodic diffeomorphism of
period equal to the order $k$ of the singularity. This is consistent with the
fact that the local fundamental group of the orbifold $L_p$ at $p$ is the cyclic
group of order $k$.

\vskip 2truecm\par
{\bf 4. Nef models}
\vglue 1truecm\par

As the title suggests, one of the purposes of [MQ1] is to develop a ``Mori
theory'' for foliations on algebraic surfaces (or even on higher dimensional
varieties). Therefore, as in the ordinary Mori theory (see e.g. [M-P]), one of
the first steps should be the construction of birational models of foliations
whose canonical bundles have suitable semipositivity properties.
\par
In what follows, $X$ will denote a complex projective surface with (at most)
cyclic quotient singularities. Recall that a ${\bf Q}$-bundle $L$ (or a ${\bf
Q}$-divisor $D$) on $X$ is called {\bf pseudoeffective} if $L\cdot H\ge 0$ for
every ample divisor $H$ on $X$. This is equivalent to say that $L\cdot C\ge 0$
for every irreducible curve $C\subset X$ whose selfintersection $C\cdot C$ is
nonnegative. If the same inequality holds for every irreducible curve $C$,
regardless its selfintersection, then $L$ is said to be {\bf nef} (numerically
eventually free). 
A useful result of Zariski and Fujita [Fuj] [Sak] states that any
pseudoeffective ${\bf Q}$-bundle $L$ has a {\it Zariski decomposition}: it can
be uniquely written as $$L=P+N$$ where $P$ (the positive part) is a nef ${\bf
Q}$-bundle and $N$ (the negative part) is an effective ${\bf Q}$-bundle whose
support is contractible and orthogonal to $P$ (i.e. $N={\cal O}_X(\sum_{j=1}^n
a_jC_j)$, with $a_j\in{\bf Q}^+$, $C_j$ irreducible curves with 
$(C_i\cdot C_j)_{i,j}$ negative definite, and $P\cdot C_j=0$ for every $j$).
\par
Let now ${\cal F}$ be a foliation on $X$, with as usual $Sing({\cal F})\cap
Sing(X)=\emptyset$, and let $K_{\cal F}\in Pic(X)\otimes{\bf Q}$ be its
canonical bundle. The following result has been proved (in a much more general
context) by Miyaoka [Miy] and Shepherd-Barron [ShB] using positive
characteristic arguments (see also [M-P, Lecture III]). 
We shall reproduce here a
characteristic zero proof due to Bogomolov and McQuillan [B-M], which reduces
the statement to a classical result of Arakelov on fibrations by curves [Ara]
[Szp] [BPV, pages 107-110].
\par
\underbar{\bf Theorem 1} [Miy] [ShB] [B-M]. {\it Let $X$ be a complex projective
surface with at most cyclic quotient singularities, and let ${\cal F}$ be a
foliation on $X$ with reduced singularities. Then the following two statements
are equivalent:
\par\noindent
i) ${\cal F}$ is a rational fibration, i.e. a fibration whose generic fibre is a
rational curve ${\bf C}P^1$;
\par\noindent
ii) $K_{\cal F}$ is not pseudoeffective.}
\par
{\it Proof}.
\par
i) $\Rightarrow$ ii). This is obvious: a generic fibre $C$ satisfies $C\cdot
C=0$ and $K_{\cal F}\cdot C=-\chi (C)=-2<0$, thus $K_{\cal F}$ cannot be
pseudoeffective.
\par
ii) $\Rightarrow$ i). If $K_{\cal F}$ is not pseudoeffective then we can find an
irreducible curve $C\subset X$ with $C\cdot C>0$ and $K_{\cal F}\cdot C <0$. We
may also suppose that $C$ is smooth, disjoint from $Sing({\cal F})\cup Sing(X)$,
and that it is not ${\cal F}$-invariant ($C$ is the zero set of a generic
section of a very ample bundle).
\par
On the projective threefold $Y=X\times C$ let us consider the 1-dimensional
foliation ${\cal G}$ which is tangent to each stratum $X\times\{ c\}$, $c\in C$,
and which coincides there with ${\cal F}$. Let $D\subset Y$ be the diagonal
curve $\{ (c,c)\vert c\in C\}\subset X\times C$. On a neighbourhood $U$ of $D$ 
we can construct a smooth complex surface $Z\subset U\subset Y$ by glueing 
together the local leaves of ${\cal G}$ through points of $D$: note that at each
$(c,c)\in D$ the foliation ${\cal G}$ is nonsingular, and its local leaf through
$(c,c)$ is a disk {\it not} tangent to $D$ (even if the local leaf of ${\cal F}$
through $c$ may be tangent to $C$). 
\par
This surface $Z$ contains $D$, and we can compute the selfintersection of $D$ in
$Z$. Because ${\cal G}$ is tangent to $Z$ and everywhere transverse to $D$ in
$Z$, we have $(D\cdot D)_Z = T_{\cal G}\cdot D$. But $T_{\cal G}$ is the 
pull-back of $T_{\cal F}$ by the projection $Y\to X$, and $D$ projects to $C$
under the same projection, thus we have
$$(D\cdot D)_Z = -K_{\cal F}\cdot C >0$$
and in particular $Z$ is a so-called {\it pseudoconcave} surface [And]. By a
theorem of Andreotti [And] (which, in this special case, is reproved in [B-M] by
an easier formal argument) the Zariski-closure of $Z$ in $Y$ has the same
dimension as $Z$: there exists an {\it algebraic} surface $W\subset Y$ 
which contains $Z$.
\par
Such a surface $W$ is ${\cal G}$-invariant (for $Z$ is), and the projection
$W\buildrel\pi\over\to C$ to the second factor of $X\times C$ is a fibration
which coincides with ${\cal G}\vert_W$ (for this is true on $Z$). The curve
$D\subset W$ is a section of this fibration, and the positivity of its
selfintersection implies, by Arakelov's theorem [Ara] [Szp], that $\pi$ is a
rational fibration. By projecting from $W$ to $X$, we see that ${\cal F}$ is 
a foliation by rational curves, and indeed a rational fibration for its
singularities are reduced (each $p\in Sing({\cal F})$ belongs to at most 2
${\cal F}$-invariant curves, thus a generic rational fibre of $\pi$ projects to
a rational leaf of ${\cal F}$, free of singularities).
\par 
$\triangle$
\par
Suppose now that the reduced foliation ${\cal F}$ is not a rational fibration,
and let us analyse the failure of nefness of its canonical bundle $K_{\cal F}$.
Thus, take an irreducible curve $C\subset X$ with $K_{\cal F}\cdot C <0$.
Because $K_{\cal F}$ is pseudoeffective, we forcely have $C\cdot C <0$, and this
implies that $C$ is ${\cal F}$-invariant: otherwise we should have $(K_{\cal
F}+C)\cdot C\ge 0$, which is not the case. We can therefore apply the formula
$$K_{\cal F}\cdot C = -\chi_{orb}(C)+Z({\cal F},C) .$$
From $K_{\cal F}\cdot C<0$ we obtain $\chi_{orb}(C) > Z({\cal F},C)\ge 0$ ($C$
has at most normal crossing singularities, because ${\cal F}$ is reduced). In
particular, $C$ is a smooth rational curve, possibly passing through $Sing(X)$.
Observe now that $C$ must intersect $Sing({\cal F})$, otherwise by Camacho--Sad
formula we would deduce $C\cdot C =0$. Hence $Z({\cal F},C)\ge 1$, and
consequently $\chi_{orb}(C) >1$. On the other side, $\chi_{orb}(C)\le 2$ and so
$Z({\cal F},C)<2$.
This leaves open only the following possibility:
$C$ is a smooth rational curve which contains at most 1 point $q$ from $Sing(X)$
(if $q$ has order $k$ then $\chi_{orb}(C)=1+{1\over k}$), and exactly 1 point
$p$ from $Sing({\cal F})$, with $Z({\cal F},C,p)=1$.
\par
The holonomy of ${\cal F}$ along a loop in $C$ around $p$ is of course equal to
the holonomy along a loop around $q$, hence it is $k$-periodic ($k=1$ if $C\cap
Sing(X)=\emptyset$). Thus ${\cal F}$ around $p$ is generated by a vector field
of the type $kz{\partial\over\partial z}-lw{\partial\over\partial w}$, with
$C=\{ w=0\}$, $l\in{\bf N}^+$ prime to $k$. By Camacho--Sad formula $l$ can be
computed from the selfintersection of $C$: we have $C\cdot C =-{l\over k}$.

\centerline{\hbox{\psfig{figure=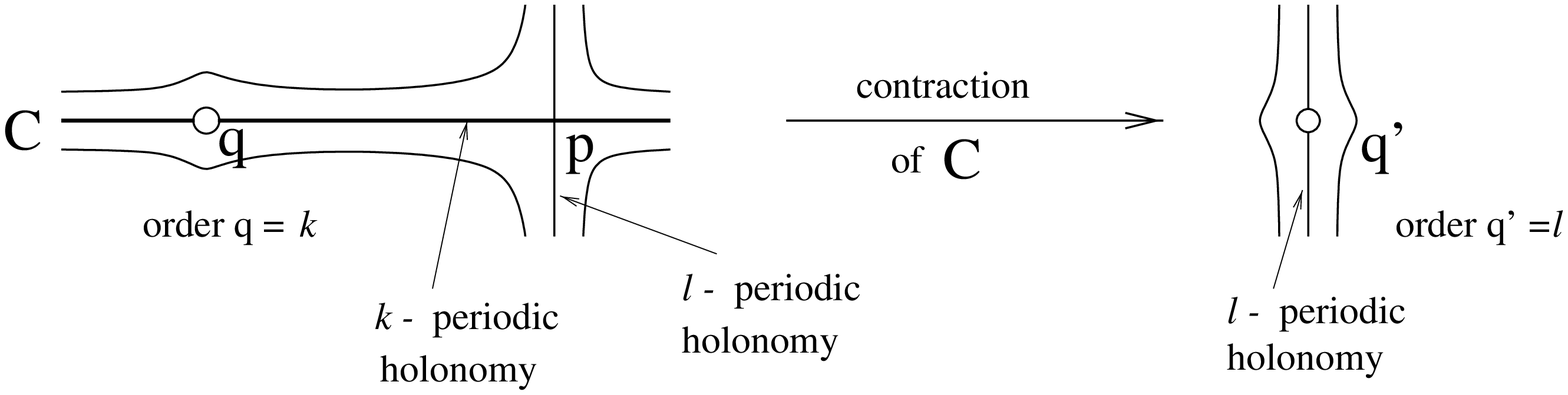,height=4truecm}}}

We can contract this curve $C$ to a point, obtaining in this way a new surface
$X'$ with a point $q'$ arising from $C$. It turns out that $q'$ is still a
cyclic quotient singularity, of order $l$ (a regular point if $l=1$), and $X'$
is still a projective surface [Fri, pages 75-76]. 
Moreover, the foliation ${\cal F}$ induces
a foliation ${\cal F}'$ on $X'$ which is nonsingular at $q'$ (and hence reduced,
as ${\cal F}$ is).
\par
In other words, by contracting a curve over which the canonical bundle of the
foliation has negative degree we obtain a new surface and a new foliation which
still satisfy our standing assumptions on $Sing(X)$ and $Sing({\cal F})$. By
iterating this procedure a finite number of steps we finally obtain the 
following result of McQuillan [MQ1] [Br1].
\par
\underbar{\bf Corollary 1}. {\it Let $X$ be a complex projective
surface with at most cyclic quotient singularities, and let ${\cal F}$ be a
foliation on $X$ with reduced singularities. Suppose that ${\cal F}$ is not a
rational fibration. Then there exists a birational morphism $(X,{\cal F})
\rightarrow (X',{\cal F}')$ such that:
\par\noindent
i) $X'$ is still projective with at most cyclic quotient singularities, and
${\cal F}'$ is still reduced;
\par\noindent
ii) the canonical bundle of ${\cal F}'$ is nef.}
\par
It is here that the presence of cyclic quotient singularities becomes
unavoidable: even if $X$ is smooth, it may happen that $X'$ is singular. Of
course, the contraction $X\to X'$ is nothing but the contraction of the support
of the negative part of the Zariski decomposition of $K_{\cal F}$ 
[MQ1, \S III.2] [Br1, page 113], thus the Corollary can be seen as a 
description of that negative part.
\par
We shall say that a foliation ${\cal F}$ is a {\bf nef foliation} if $K_{\cal
F}$ is nef, so that the previous results (Seidenberg's resolution, Miyaoka's
theorem, McQuillan's contraction) say that {\it any foliation, which is
not birational to a rational fibration, has a birational model which is reduced
and nef}. Such a model, in general, is not unique: for instance, the blow-up of
a reduced nef foliation at a singular point is still reduced and nef (whereas
nefness is lost if we blow up a regular point). In order to get uniqueness we
need, at least, to contract more curves, and to obtain the so-called {\it
minimal models}; but even after these additional contractions uniqueness may
fail, and we refer to [MQ1] [Br1] [Br3] for a more complete discussion of this
point and related ones.
\par
From now on, we shall concentrate our attention to {\it reduced nef foliations}.
Given such a foliation ${\cal F}$, with canonical bundle $K_{\cal F}\in
Pic(X)\otimes{\bf Q}$, we define the {\bf numerical Kodaira dimension} $\nu
({\cal F})$ of ${\cal F}$ as the numerical Kodaira dimension of $K_{\cal F}$,
that is:
\par\noindent
$\bullet$ $\nu ({\cal F})=2$ if $K_{\cal F}\cdot K_{\cal F} >0$;
\par\noindent
$\bullet$ $\nu ({\cal F})=1$ if $K_{\cal F}\cdot K_{\cal F} =0$ but 
$K_{\cal F}$ is not numerically trivial (there exists $C$ such that 
$K_{\cal F}\cdot C >0$);
\par\noindent
$\bullet$ $\nu ({\cal F})=0$ if $K_{\cal F}$ is numerically trivial 
($K_{\cal F}\cdot C=0$ for every $C$).
\par\noindent
(Because $K_{\cal F}$ is nef, we have $K_{\cal F}\cdot K_{\cal F}\ge 0$).
\par
If $\nu ({\cal F})=2$ then ${\cal F}$ is said to be of {\bf general type}. In
some sense, most foliations belong to this class. Our ultimate task will be the
classification of foliations which are {\it not} of general type.

\vskip 2truecm\par
{\bf 5. Numerically trivial foliations}
\vglue 1truecm\par

A reduced nef foliation ${\cal F}$ on a projective surface $X$ is said to be
{\bf numerically trivial} if $\nu ({\cal F})=0$, that is its canonical bundle 
$K_{\cal F}$ has zero degree on every curve $C\subset X$. 
The classification of these foliations is rather easy: not surprisingly, 
they are related to foliations whose
canonical bundle is {\it holomorphically} trivial, i.e. foliations which are
globally generated by a single global holomorphic vector field with isolated 
zeroes.
\par
\underbar{\bf Theorem 2} [MQ1] [Br1]. {\it Let ${\cal F}$ be a reduced nef 
foliation with $$\nu ({\cal F})=0$$ on a projective surface $X$. 
Then there exists a
finite regular covering $Y\buildrel\pi\over\rightarrow X$ such that the lifted
foliation ${\cal G}=\pi^*({\cal F})$ has a holomorphically trivial canonical
bundle $K_{\cal G}$.}
\par
(Note: $Y$ is smooth, and {\it regular covering} has to be understood in
orbifold's sense, i.e. on a neighbourhood of a point $q\in Y$ sent to $p\in
Sing(X)$ the map $\pi$ looks like ${\bf B}^2\to {\bf B}^2/\Gamma$).
\par
{\it Proof.}
\par
The main step consists in showing that $K_{\cal F}$ is a torsion ${\bf
Q}$-bundle (i.e., some positive power of $K_{\cal F}$ is effective),
then the conclusion will follow by the standard covering trick [MQ1,
\S IV.3] [Br1, pages 110-112]. 
\par
If $H^1(X,{\cal O}_X)=0$ then there is nothing to prove: on such a surface, any
numerically trivial ${\bf Q}$-bundle is forcely a torsion ${\bf Q}$-bundle,
because the Chern class map $H^1(X,{\cal O}_X^*)\rightarrow H^2(X,{\bf Z})$ is
injective and the numerical triviality is equivalent to the vanishing of the
rational Chern class in $H^2(X,{\bf Q})$.
\par
Thus let us suppose $H^1(X,{\cal O}_X)\not= 0$, whence $H^0(X,\Omega_X^1)\not=
0$ by Hodge symmetry. If some global 1-form does not vanish identically
when restricted to the leaves of ${\cal F}$, then this restriction defines a
nontrivial section of $K_{\cal F}$, and we are done (in this case $K_{\cal F}$
is already trivial and $Y=X$). If every global 1-form vanishes on the leaves of
${\cal F}$, then we necessarily are in the following situation: the Albanese map
$X\buildrel alb\over\rightarrow Alb(X)$ has a 1-dimensional image $B$, and
${\cal F}$ coincides with the fibration $X\buildrel alb\over\rightarrow B$. The
generic fibres are elliptic curves and the fibration is isotrivial, for $K_{\cal
F}$ is numerically trivial (Arakelov's theorem [Ara] [Ser]). Then it is still a
consequence of results of Arakelov that $K_{\cal F}$ is a torsion 
${\bf Q}$-bundle [MQ1, page 65] [Br1, page 111].
\par
$\triangle$
\par
Recall now that global holomorphic vector fields on projective surfaces are well
understood, see for instance [Br1, Chapter 6]. If $v$ is such a vector field,
with isolated singularities, on a surface $Y$, then we have one of the following
four possibilities:
\par\noindent
1) $Y$ has an isotrivial elliptic fibration, all of whose fibres have smooth
reduction, and $v$ is tangent to the fibres, nowhere vanishing;
\par\noindent
2) $Y$ is a torus ${\bf C}^2/G$, and $v$ is the quotient of a constant vector
field on ${\bf C}^2$;
\par\noindent
3) $Y$ is a ${\bf C}P^1$-bundle over an elliptic curve $E$, and $v$ is
transverse to the fibres and projects on $E$ to a constant vector field;
\par\noindent
4) $Y$ is a rational surface, and up to a birational map we have $Y={\bf C}P^1
\times{\bf C}P^1$ and $v=v_1\oplus v_2$, $v_1$ and $v_2$ being holomorphic
vector fields on ${\bf C}P^1$.
\par
These facts and Theorem 2 give a rather exhaustive description of reduced nef
foliations with vanishing numerical Kodaira dimension. Let us just observe the
following corollary to Theorem 2: if $\nu ({\cal F})=0$ then every leaf of
${\cal F}$ is {\it parabolic}, i.e. uniformized by ${\bf C}$. We shall see later
another way to prove such a corollary, using the Poincar\'e metric on the
leaves. Be careful, however, that there are foliations with $\nu ({\cal F})=1$
which also have all their leaves parabolic: see for instance [Br4] for some
natural examples. This is quite unpleasant.

\vskip 2truecm\par
{\bf 6. Kodaira dimension}
\vglue 1truecm\par

The {\bf Kodaira dimension} $kod({\cal F})$ of a reduced nef foliation ${\cal
F}$ on a projective surface $X$ is defined as the Kodaira-Iitaka dimension of
its canonical bundle $K_{\cal F}\in Pic(X)\otimes{\bf Q}$, that is
$$kod({\cal F})=\limsup_{n\to +\infty}{\log\dim H^0(X,K_{\cal F}^{\otimes n})
\over \log n} \in \{ -\infty ,0,1,2\}.$$
Standard arguments give $kod({\cal F})\le\nu ({\cal F})$, and moreover
$kod({\cal F})=2$ $\Longleftrightarrow$ $\nu ({\cal F})=2$. In the proof of
Theorem 2 we showed $\nu ({\cal F})=0$ $\Longrightarrow$ $kod({\cal F})=0$. Here
we shall prove the converse.
\par
\underbar{\bf Theorem 3} [MQ1]. {\it Let ${\cal F}$ be a reduced nef foliation
with $\nu ({\cal F})=1$, on a projective surface $X$. Then $kod({\cal F})$ is
either $-\infty$ or $1$.}
\par
{\it Proof}.
\par
We shall give a proof following [Br1, pages 119-126], which is slightly 
different from [MQ1].
\par
Suppose $kod({\cal F})\ge 0$, i.e. for some positive integer $n$ the
canonical bundle $K_{\cal F}^{\otimes n}$ has a nontrivial global section $s$
vanishing on an effective divisor $D$. We have $D\not= \emptyset$ and $D\cdot
D=0$, for $\nu({\cal F})=1$. For each irreducible component $D_j$ of $Supp(D)$ 
we have $D\cdot D_j=0$ and $D_j^2\le 0$, for $D$ is nef. 
\par
{\bf 1)}. If some component $D_j$ is not ${\cal F}$-invariant then, from
$0=K_{\cal F}\cdot D_j=$ $-D_j^2+tang({\cal F},D_j)$, we find that $D_j$ is
everywhere transverse to ${\cal F}$ and $D_j^2=0$, in particular $D_j$ is
disjoint from the other components of $D$. The section $s$ of $K_{\cal
F}^{\otimes n}$ induces on each local leaf $L_p$ of ${\cal F}$ through $p\in
D_j$ a section of $\Omega^1(L_p)^{\otimes n}$, vanishing only at $p$. By
integrating this section on each local leaf [Br1, page 121] we obtain a
holomorphic function $f$ on a neighbourhood of $D_j$, vanishing on $D_j$ and
only there. The level sets of $f$ define a fibration which contains $D_j$ as a
(multiple) fibre, and standard compactness arguments show that such a fibration
extends to the full $X$. Thus there exists a fibration
$X\buildrel\pi\over\rightarrow B$ such that $D_j=\pi^*(Z)$ where $Z$ is a
positive ${\bf Q}$-divisor on $B$, supported on a single point (the image of
$D_j$ by $\pi$). It follows that $kod({\cal O}_X(D_j))=kod({\cal O}_B(Z))=1$,
whence $kod({\cal F})=kod({\cal O}_X(D))=1$.
\par
{\bf 2)}. From now on we shall therefore assume that each irreducible component
$D_j$ is ${\cal F}$-invariant. From $0=K_{\cal F}\cdot D_j=$
$-\chi_{orb}(D_j)+Z({\cal F},D_j)$ we thus find
$$\chi_{orb}(D_j)=Z({\cal F},D_j).$$
Remark that each point of $D_j\cap D_k$, $k\not= j$, is a singularity of ${\cal
F}$ and gives a positive contribution to $Z({\cal F},D_j)$. Hence, setting
$D_{red}=\cup_iD_i$ we have
$$Z({\cal F},D_j)\ge D_{red}\cdot D_j - D_j^2$$
and consequently
$$K_X\cdot D_j = -D_j^2 - \chi_{orb}(D_j)\le -D_{red}\cdot D_j$$
$$K_X\cdot D \le -D_{red}\cdot D =0 .$$
If $K_X\cdot D < 0$ then by Riemann--Roch formula we deduce $kod({\cal O}_X(D))
=1$, as desired. Hence we shall assume from now on that
$$K_X\cdot D =0 .$$
Thus the previous inequalities become equalities, that is $Z({\cal F},D_j)=
D_{red}\cdot D_j-D_j^2$ for every $j$. This means that ${\cal F}$ is singular on
$D_j$ only in correspondence of the intersections with the other components, or
in correspondence of its selfintersections. Moreover all these singularities are
nondegenerate ($Z({\cal F},D_j,p)=Z({\cal F},D_k,p)=1$ if $p\in D_j\cap D_k$).
\par
Let now $D_0$ be a connected component of $Supp(D)$, not necessarily
irreducible. By the previous considerations, we have $Sing({\cal F})\cap D_0=$
nodal points of $D_0$. From this and $\chi_{orb}(D_j)=Z({\cal F},D_j)$ for every
$j$ it is then easy to find the structure of $D_0$ (see the discussion before
Corollary 1). Here is a list of all the possibilities (compare with Kodaira's
table of elliptic fibres [BPV, page 150] [Br1, page 67]: case $(a)$ corresponds 
to $I_0$, case $(b)$ to $I_b$, $b\ge 1$, case $(c)$ to $II$, $II^*$, $III$,
$III^*$, $IV$ or $IV^*$, case $(d)$ to $I_0^*$, case $(e)$ to $I_b^*$, $b\ge
1$):
\par\noindent
$(a)$ $D_0$ is a smooth elliptic curve, disjoint from $Sing({\cal F})$ and
$Sing(X)$;
\par\noindent
$(b)$ $D_0$ is a cycle of smooth rational curves, or a rational curve with a
node; $D_0\cap Sing(X)=\emptyset$ and $D_0\cap Sing({\cal F})=$ nodes of $D_0$;
\par\noindent
$(c)$ $D_0$ is a rational curve, disjoint from $Sing({\cal F})$ and with
$D_0\cap Sing(X)=\{ q_1,q_2,q_3\}$; the orders $k_j$ of $q_j$ satisfy ${1\over
k_1}+{1\over k_2}+{1\over k_3}=1$;
\par\noindent
$(d)$ $D_0$ is a rational curve, disjoint from $Sing({\cal F})$ and with
$D_0\cap Sing(X)=\{ q_1,q_2,q_3,q_4\}$; each $q_j$ has order 2;
\par\noindent 
$(e)$ $D_0$ is a chain of rational curves, with $D_0\cap Sing({\cal F})=$ nodes
of $D_0$; each extreme curve in $D_0$ contains two points from $Sing(X)$, both
of order 2; each interior curve in $D_0$ is disjoint from $Sing(X)$.

\centerline{\hbox{\psfig{figure=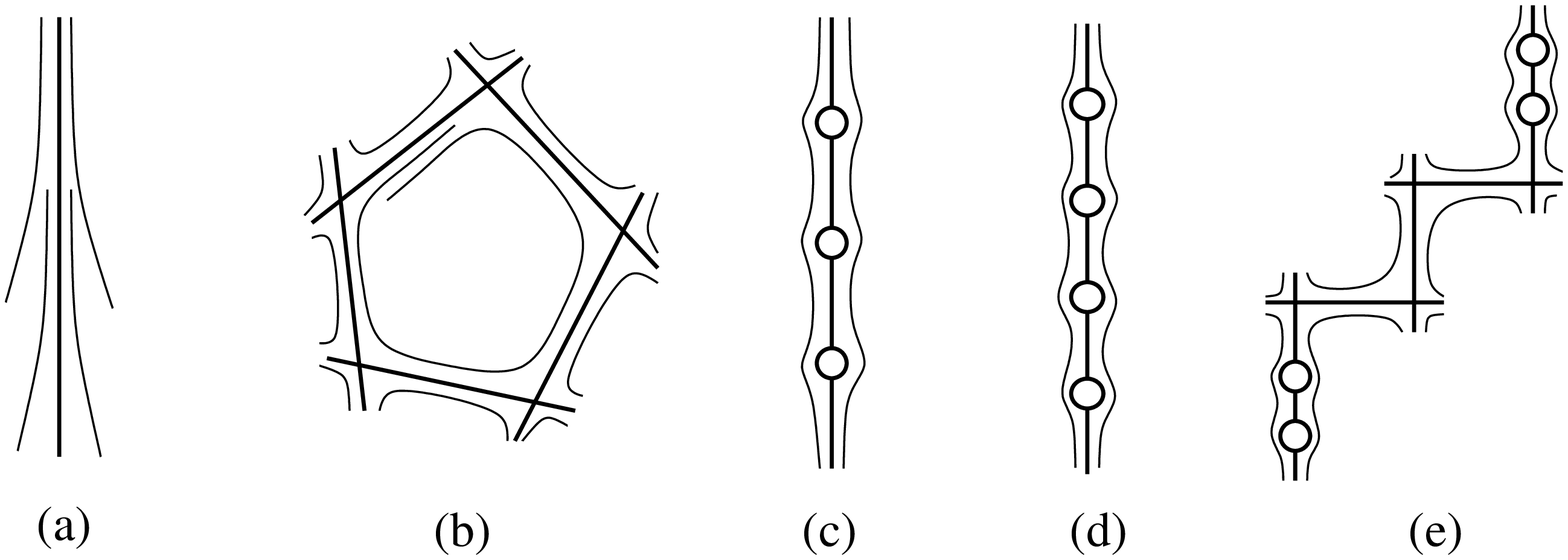,height=5truecm}}}

In case $(b)$ there could be a curve of selfintersection $-1$ among the rational
curves of the cycle, but then we may contract it and we obtain a similar and
simpler situation. Thus, we may suppose that all the curves in $(b)$ have
selfintersection $-2$ (recall that $D_0$ is a connected component of $Supp(D)$
and $D$ is nef with zero selfintersection). Same thing in case $(e)$: we may
suppose that the interior curves have selfintersection $-2$, the extreme ones
$-1$. In both cases the singularities of ${\cal F}$ on $D_0$ can be computed via
Camacho--Sad formula: the $CS$ indices must be all equal to $-1$, hence the
singularities are all of the type $z{\partial\over\partial z} -
w(1+...){\partial\over\partial w}$.
\par
Remark also that in all the cases we have $D=mD_0$ ($+$ other components), i.e.
all the irreducible components of $D_0$ appear in $D$ with the same multiplicity
$m$ (this follows again from $D$ nef, $D^2=0$). Hence $D_0^2=0$.
\par
As in the first part of the proof (with one irreducible component not ${\cal
F}$-invariant), we need just to prove that $D_0$ is a (multiple) fibre of some
(elliptic) fibration $X\buildrel\pi\over\rightarrow B$.
\par
{\bf 3)}. If $H^1(X,{\cal O}_X)\not= 0$ we consider, as in Theorem 2, the
Albanese map $X\buildrel alb\over\rightarrow Alb(X)$. Because $D_0^2=0$, we have
the alternative: either $alb$ is a fibration over a curve and $D_0$ is a fibre,
or $D_0$ is not contracted by $alb$ and so there exists a 1-form
$\omega\in\Omega^1(X)$ with $\omega\vert_{D_0}\not\equiv 0$. In the former case
we have obviously finished (with $\pi =alb$). In the latter case
$\omega\vert_{\cal F}$ defines a section of $K_{\cal F}$ not vanishing on $D_0$,
which easily gives $kod({\cal F})=1$. Therefore we shall assume from now on that
$$H^1(X,{\cal O}_X)=0.$$
We may also assume that $Supp(D)$ is connected, i.e. $$D=mD_0.$$
Indeed, if there exists a second connected component $D_0'$ then, by Hodge index
theorem [BPV] and $H^1(X,{\cal O}_X)=0$, we have ${\cal O}_X(aD_0)={\cal
O_X}(bD_0')$ for some positive $a$ and $b$, whence $kod({\cal O}_X(D))=1$.
\par
{\bf 4)}. Let $l$ be the minimal positive integer such that $lD_0$ is a Cartier
divisor: $l=1$ in cases $(a)$ and $(b)$, $l=2$ in cases $(d)$ and $(e)$,
$l=m.c.m.(k_1,k_2,k_3)\in\{ 3,4,6\}$ in case $(c)$. By Riemann--Roch formula 
(and $D_0^2=K_X\cdot D_0=0$, $\chi (X,{\cal O}_X)\ge 1$) we have, for $n\gg 0$,
$$h^0(X,{\cal O}_X(nlD_0))\ge h^1(X,{\cal O}_X(nlD_0)) +1.$$
If $h^1(X,{\cal O}_X(nlD_0))\ge 1$ for some $n$ then $h^0(X,{\cal O}_X(nlD_0))
\ge 2$ and so $kod({\cal O}_X(D))=1$. If $h^1(X,{\cal O}_X(nlD_0))=0$ for every
$n$ then we consider the exact sequences
$$0\rightarrow {\cal O}_X((n-1)lD_0)\rightarrow {\cal O}_X(nlD_0)\rightarrow
{\cal O}_X(nlD_0)\vert_{lD_0}\rightarrow 0$$
from which we see that
$$h^0(X,{\cal O}_X(nlD_0))\ge h^0(lD_0,{\cal O}_X(nlD_0)\vert_{lD_0}) +1.$$
Hence we are reduced to prove that ${\cal O}_X(nlD_0)$ {\it has a nontrivial
section over} $lD_0$, for some  (large) positive $n$. This is the same as to
prove that ${\cal O}_X(lD_0)\vert_{lD_0}$ is a torsion line bundle.
\par
{\bf 5)}. Recall that $D=mD_0$ is the zero set of a section $s$ of $K_{\cal
F}^{\otimes n}$ for some $n>0$. However, by the branched covering trick [MQ1,
page 83] we may assume that $n=1$, i.e. $s$ is a section of $K_{\cal F}$.
Because $K_{\cal F}$ is not locally free at points of $Sing(X)$, the section $s$
necessarily vanishes on $Sing(X)$, and so $Sing(X)\subset D_0$.
\par
Let us consider now the logarithmic conormal bundle $N_{\cal F}^*\otimes {\cal
O}_X(D_0)$, and observe that it is a genuine line bundle and not simply a ${\bf
Q}$-bundle (the logarithmic 1-form ${dz\over z}$ is invariant by $(z,w)\mapsto
(e^{{2\pi i\over k}}z,e^{{2\pi i\over k}h}w)$). By Riemann--Roch formula (and
$D_0^2=K_X\cdot D_0=0$, $N_{\cal F}^*=K_X-mD_0$) we have, for $n\gg 0$,
$$h^0(X,N_{\cal F}^*\otimes{\cal O}_X(D_0)\otimes{\cal O}_X(nlD_0))\ge 1$$
that is $N_{\cal F}^*\otimes{\cal O}_X(D_0)\otimes{\cal O}_X(nlD_0) =
{\cal O}_X(rD_0+R)$ for some effective Cartier divisor $rD_0+R$, with 
$D_0\not\subset Supp(R)$. 
In fact, we have $D_0\cap Supp(R)=\emptyset$, for otherwise we would
have $R\cdot D_0 >0$ which contradicts $N_{\cal F}^*\cdot D_0=0$, $D_0^2=0$.
Because $rD_0+R$ is Cartier, we therefore have $r=$ multiple of $l$. Thus we
finally obtain, for some $q\in {\bf Z}$:
$$N_{\cal F}^*\otimes{\cal O}_X(D_0) = {\cal O}_X(qlD_0+R)$$
where $R$ is a Cartier divisor disjoint from $D_0$. In other words, ${\cal F}$
is globally defined by a meromorphic 1-form $\omega$ with poles along $D_0$ of
order $1-ql$.
\par
The case $q>0$ is excluded: $\omega$ would be holomorphic, contradicting
$H^0(X,\Omega_X^1)=$ $H^1(X,{\cal O}_X)=0$. Also the case $q=0$ is excluded:
$\omega$ would be a logarithmic 1-form with nontrivial residue $\sum\lambda_j
D_j$, but this is impossible because the residue of a logarithmic 1-form must be
cohomologous to zero (in the formalism of currents, $Res(\omega )={1\over 2\pi
i}\bar\partial\omega$), which is certainly not the case of $\sum\lambda_j D_j$
unless $\lambda_j=0$ for every $j$. Thus we have $q<0$.
\par
{\bf 6)}. If $l=1$ (cases $(a)$ and $(b)$) the proof is easily achieved. Indeed,
we claim that $N_{\cal F}^*\otimes{\cal O}_X(D_0)\vert_{D_0}$ is trivial. To see
this, observe that each local section of $N_{\cal F}^*\otimes{\cal O}_X(D_0)$
has a residue on $D_0$. If $p$ is a smooth point of $D_0$ then the
residue is, around $p$, a holomorphic function on $D_0$. If $p$ is a nodal point
of $D_0$ then the residue is, around $p$, a pair of holomorphic functions on the
two local branches of $D_0$ through $p$; however, using the fact that the
singularities of ${\cal F}$ are of the type $zdw+w(1+...)dz=0$, we see that
those two functions coincide at $p$, and so they can be identified with a single
function on $D_0$. This means that we have a map 
$$N_{\cal F}^*\otimes{\cal O}_X(D_0)\buildrel Res\over\longrightarrow 
{\cal O}_{D_0}$$
and this map becomes an isomorphism when we restrict to $D_0$.
\par
From the triviality of $N_{\cal F}^*\otimes{\cal O}_X(D_0)\vert_{D_0}$ and the
previous step it now follows that 
${\cal O}_X(qD_0+R)\vert_{D_0}={\cal O}_X(qD_0)\vert_{D_0}$ is trivial too 
and has a nonvanishing section. 
\par
If $l>1$ we need a more deep analysis of the structure of ${\cal F}$ around 
$D_0$, which takes into account higher order terms given by the holonomy of the
foliation.
\par
{\bf 7)}. If $l>1$ (cases $(c)$, $(d)$ and $(e)$) we may do, on a neighbourhood
$U$ of $D_0$, the semistable reduction [Br1, page 125] [BPV, page 155]: there
exists a regular covering of order $l$, $V\buildrel r\over\rightarrow U$ ($V$
smooth, $r$ regular in orbifold's sense), such that $E_0=r^{-1}(D_0)$ is a curve
of type $(a)$ or $(b)$ and ${\cal G}=r^*({\cal F})$ is tangent to $E_0$. By
pulling back the meromorphic 1-form $\omega$ we see that
$$N_{\cal G}^*\otimes{\cal O}_V(E_0) = {\cal O}_V(qlE_0).$$
\par
Let us firstly consider the case in which $D_0$ is of type $(c)$ or $(d)$, so
that $E_0$ is of type $(a)$. The holonomy group $Hol({\cal F})\subset 
Diff({\bf C},0)$ of ${\cal F}$ along $D_0$ is a solvable group, generated by the
finite order holonomies around points in $Sing(X)$. If this group is finite then
${\cal F}$, around $D_0$ and hence everywhere, is a fibration, and the proof is
finished. If this group is not finite then [C-M] [L-M] its commutator is
an abelian infinite group which embeds in a flow of a germ of vector field of
the type $z^{p+1}{\partial\over\partial z}$, where $p\in{\bf N}^+\setminus l{\bf
N}^+$. This commutator is in fact the holonomy group of ${\cal G}$ along $E_0$,
and the previous property can be used to construct, around $E_0$, a (closed)
meromorphic 1-form which defines ${\cal G}$ and has a pole of order $p+1$ along
$E_0$ [Br1, page 126] [Pau]. This means, in particular, that
$$N_{\cal G}^*\otimes{\cal O}_V(E_0) = {\cal O}_V(-pE_0).$$
Hence
$${\cal O}_V(qlE_0)={\cal O}_V(-pE_0)$$
and from $p\not= -ql$ we infer that ${\cal O}_V(E_0)$ is a torsion line bundle.
Then the same holds for ${\cal O}_U(lD_0)$, so that ${\cal
O}_X(nlD_0)\vert_{lD_0}$ is trivial and has a nonvanishing section for some
$n>0$.
\par
Let us finally consider the case in which $D_0$ is of the type $(e)$, so that
$E_0$ is of type $(b)$. Recall that the singularities of ${\cal F}$ on $D_0$ are
of the type $zdw+w(1+...)dz=0$. If one singularity has a first integral, i.e.
trivial holonomy, then it is easy to see that the same holds for any other
singularity and ${\cal F}$, around $D_0$ and hence everywhere, is a fibration.
Thus we are left with the case in which no singularity has a first integral.
Each extreme curve of $D_0$ has then a solvable nonabelian holonomy, and as
before we find [C-M] that the singularity of ${\cal F}$ on that extreme
curve has holonomy embeddable (at least formally) in the flow of 
$z^{p+1}{\partial\over\partial z}$ for some $p\in{\bf N}^+\setminus 2{\bf N}^+$
(in fact, this is just a simple consequence of the fact that this holonomy is
the product of two 2-periodic holonomies).
This means [MR2] that, in suitable (formal) coordinates around that point, 
${\cal F}$ is defined by the closed meromorphic 1-form 
$[{1\over (zw)^p}+\lambda -1]{dz\over z} + 
[{1\over (zw)^p}+\lambda ]{dw\over w}$, for some $\lambda\in{\bf C}$ (in fact
here $\lambda$ is zero). This property propagates to the other singularities,
with the same $p$ (and $\lambda$ integer). We now pass to the covering ${\cal
G}$, where we find the same singularities, and then [Pau] we may construct (at
least formally) a (closed) meromorphic 1-form which defines ${\cal G}$ and has 
a pole of order $p+1$ along $E_0$. As before, because $p\not= 2$ this gives 
the triviality of ${\cal O}_X(2nD_0)\vert_{2D_0}$ for some $n>0$ (when we
restrict to $2D_0$ all the formal problems disappear, of course).
\par
$\triangle$
\par
We shall see later that foliations with $\nu ({\cal F})=1$ and $kod({\cal F})=
-\infty$ indeed exist. In the next section we address to the much easier case
$kod({\cal F})=1$.

\vskip 2truecm\par
{\bf 7. Riccati and Turbulent foliations}
\vglue 1truecm\par

A foliation ${\cal F}$ on $X$ is a {\bf Riccati foliation}, resp. a {\bf
Turbulent foliation}, if there exists a fibration $\pi :X\to B$ whose generic
fibre is rational, resp. elliptic, and transverse to ${\cal F}$. Riccati and
Turbulent foliations constitute a basic piece of the classification of
nongeneral type foliations, as the next result shows.
\par
\underbar{\bf Theorem 4} [MQ1] [Men]. {\it Let ${\cal F}$ be a reduced nef 
foliation with $$\nu ({\cal F})=kod({\cal F})=1$$ on a projective surface $X$.
Then one of the following situations occurs:
\par\noindent
i) ${\cal F}$ is a Riccati foliation;
\par\noindent
ii) ${\cal F}$ is a Turbulent foliation;
\par\noindent
iii) ${\cal F}$ is a nonisotrivial elliptic fibration;
\par\noindent
iv) ${\cal F}$ is an isotrivial fibration of genus $\ge 2$.}
\par
{\it Proof.}
\par
Let  $\pi :X\rightarrow B$ be the Iitaka fibration [M-P] associated to
$K_{\cal F}$. A generic fibre $C$ of $\pi$ is a curve over which $K_{\cal F}$
has zero degree. If ${\cal F}$ coincides with the Iitaka fibration then
$0=K_{\cal F}\cdot C = -\chi (C)$, thus $C$ is elliptic. The nonisotriviality of
$\pi$ follows from Arakelov's theorem [Ara] [Ser], and we are in case iii).
\par
If ${\cal F}$ does not coincides with the Iitaka fibration, then $0=K_{\cal F}
\cdot C =$ $-C^2+tang({\cal F},C)=tang({\cal F},C)$, thus $C$ is transverse to
${\cal F}$. If $C$ is rational, resp. elliptic, then ${\cal F}$ is Riccati,
resp. Turbulent. 
\par
Suppose now that $C$ has genus $\ge 2$. The transversality of
${\cal F}$ to the generic fibres of $\pi$
shows that the fibration $\pi$ is isotrivial, and so it is a
locally trivial fibre bundle up to ramified coverings and birational maps [BPV,
Chapter III]. Standard hyperbolic arguments (finiteness of $Aut({\rm fibre})$,
Picard's theorem) imply that ${\cal F}$, after those ramified coverings and
birational maps, becomes a second fibration everywhere transverse to $\pi$. In
other words, the universal covering of the surface is ${\bf D}\times\Sigma$,
with $\Sigma = {\bf D}$ or ${\bf C}$ or ${\bf C}P^1$, and on that universal
covering the foliation is the vertical one, with leaves $\{ * \}\times\Sigma$.
It follows easily that ${\cal F}$ itself is an isotrivial fibration,
the isotriviality being given by the fibres of $\pi$. Finally, $K_{\cal F}$ 
has strictly positive degree on a fibre of ${\cal F}$, because that fibre is
transverse to the Iitaka fibration of $K_{\cal F}$. Thus ${\cal F}$ has genus
$\ge 2$, and we are in case iv).
\par
$\triangle$
\par
A partial converse to this theorem is also true: nonisotrivial elliptic
fibrations and isotrivial fibrations of genus $\ge 2$ have $\nu ({\cal F})
=kod({\cal F})=1$ [Ser], Riccati and Turbulent foliations have $\nu ({\cal F})
=kod({\cal F}) \le 1$ (see below). 
\par
Let us now rapidly describe the structure of Riccati and Turbulent foliations,
following [Br1, Chapter 4] [Br4] [MQ1, \S IV.4].
\par
Let ${\cal F}$ be a reduced Riccati foliation on $X$. Then up to a birational
morphism $X\to X'$ (contraction of curves in fibres) we have that around each
fibre $F$ of the rational fibration the foliation ${\cal F}$ has the structure
given by one of the following models:

\centerline{\hbox{\psfig{figure=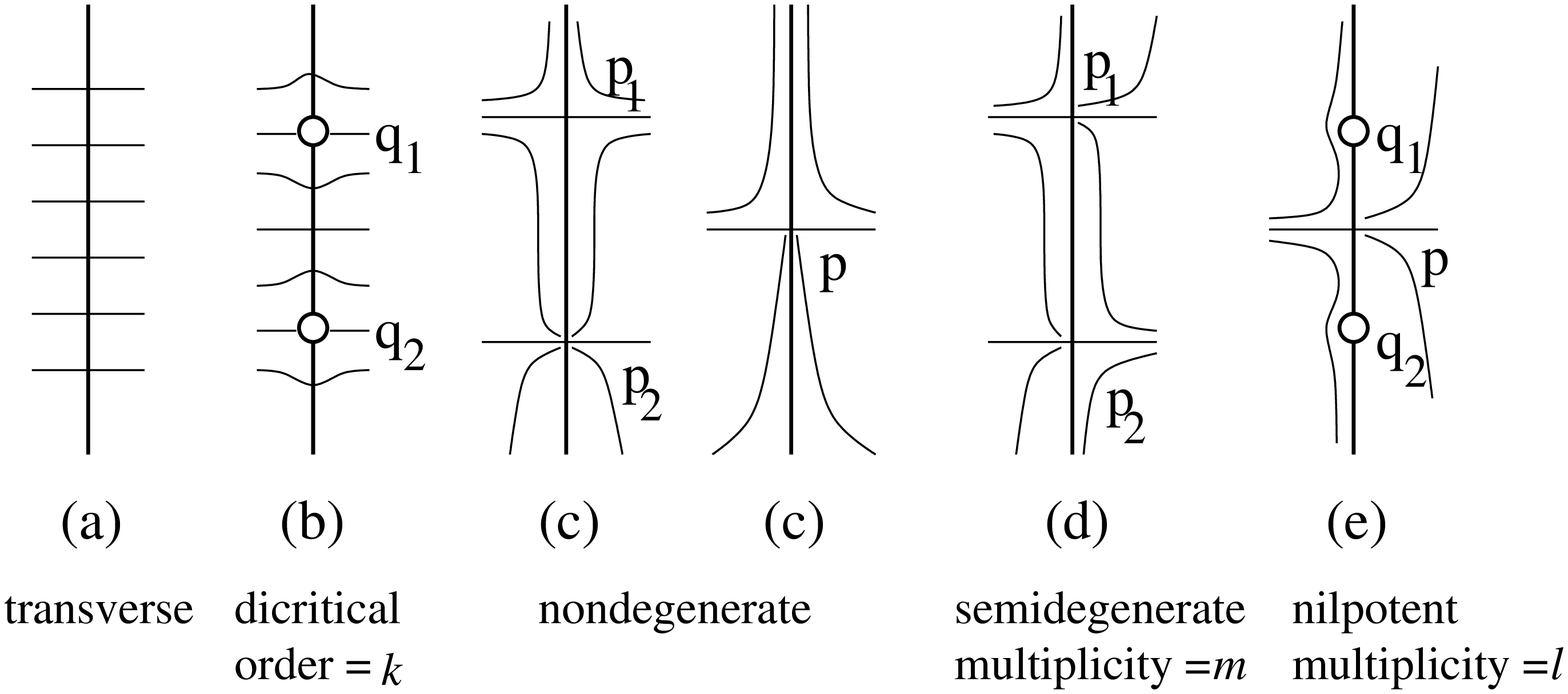,height=5truecm}}}

In $(b)$, $q_1$ and $q_2$ are quotient singularities of the same order $k$. In
$(c)$, $p_1$ and $p_2$ are nondegenerate, $p$ is a saddle-node of multiplicity 
2, with strong separatrix transverse to the fibre. In $(d)$, $p_1$ and $p_2$ are
saddle-nodes of the same multiplicity $m$, with strong separatrices inside the
fibre. In $(e)$, $q_1$ and $q_2$ are quotient singularities of order 2, $p$ is a
saddle-node of multiplicity $l$, with strong separatrix inside the fibre.
\par
By pulling back sections of $K_B$ under the fibration $X\buildrel\pi\over\to B$,
we obtain sections of $K_{\cal F}$. We can in this way compute $K_{\cal F}$,
which is the pull-back of its direct image $\pi_*(K_{\cal F})\in Pic(B)\otimes
{\bf Q}$. We finally obtain
$$deg (\pi_*K_{\cal F})= -\chi_{orb}(B) + \sum_{(c)}1 + \sum_{(d)}m_j +
\sum_{(e)}{l_j\over 2}$$
where the sums are over fibres of class $(c)$, $(d)$ and $(e)$, and the orbifold
structure of $B$ is given by assigning the multiplicity $k$, resp. 2, to a point
over which the fibre is of class $(b)$, resp. $(e)$. Hence
$$\chi_{orb}(B)=\chi_{top}(B)-\sum_{(b)}{k_j-1\over k_j}-\sum_{(e)}{1\over 2}.$$
We have $kod({\cal F})=1$ if $deg (\pi_*K_{\cal F})>0$, $kod({\cal F})=0$ if 
$deg (\pi_*K_{\cal F})=0$, $kod({\cal F})=-\infty$ if $deg (\pi_*K_{\cal F})<0$.
It is a good exercise to check that if $deg (\pi_*K_{\cal F})<0$ then ${\cal F}$
is a rational fibration, and if $deg (\pi_*K_{\cal F})=0$ then, up to a
covering, ${\cal F}$ is defined by a global vector field. Note also that it is
well possible that $deg (\pi_*K_{\cal F})>0$ but all the leaves of ${\cal F}$
are parabolic.
\par
If we set $B_0=B\setminus\{$points over which the fibre is of class $(c)$ $(d)$
$(e)\}$ then the foliation induces a monodromy representation $\pi_1^{orb}(B_0)
\rightarrow Aut({\bf C}P^1)$ (recall that points over which the fibre is of
class $(b)$ have a local fundamental group of order $k$). For each leaf $L$
outside the fibres of class $(c)$, $(d)$, $(e)$ the map
$L\buildrel\pi\over\rightarrow B_0$ is a regular covering (in orbifold's sense),
and so the orbifold universal covering of $L$ is equal to the one of $B_0$.
\par
Consider now the case of a reduced Turbulent foliation. Here the main remark is
that the elliptic fibration $X\buildrel\pi\over\rightarrow B$ is isotrivial, the
isotriviality being defined by the transverse foliation. Hence fibres of type
$I_b$ and $I_b^*$, $b\ge 1$, cannot appear [BPV, pages 150-159]. Up to a
birational morphism $X\to X'$ we then have one of the following models:

\centerline{\hbox{\psfig{figure=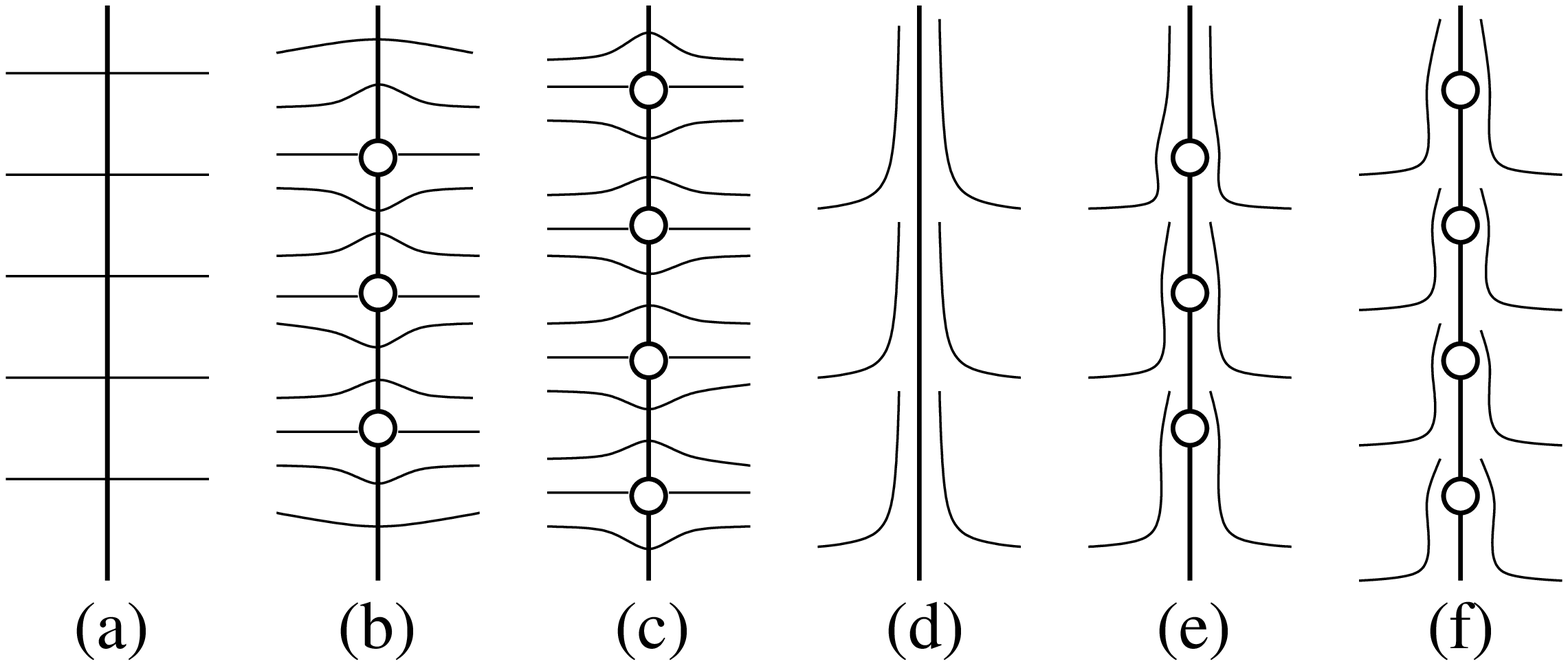,height=4truecm}}}

In cases $(a)$ and $(d)$ the fibre is smooth elliptic, and may be multiple.
In cases $(b)$ and $(e)$ the fibre is rational with three quotient singularities
of orders $k_1$, $k_2$, $k_3$, with ${1\over k_1}+{1\over k_2}+{1\over k_3}=1$;
its multiplicity is $m.c.m.(k_1,k_2,k_3)\in\{ 3,4,6\}$.
In cases $(c)$ and $(f)$ the fibre is rational with four quotient singularities
of order 2; its multiplicity is 2. In cases $(d)$, $(e)$ and
$(f)$ the foliation may be tangent to the fibre at any order. One can easily
find an explicit formula for $deg (\pi_*K_{\cal F})$, similar to the one for the
Riccati case. This allows to compute $kod({\cal F})$. We have a natural orbifold
structure on the base $B$, and a monodromy representation $\pi_1^{orb}(B_0)
\rightarrow Aut(E)$, where $E$ is an elliptic curve.
\par
By stable reduction [BPV, page 95] we may further simplify the situation: on a
regular covering $Y\to X$ the elliptic fibration becomes a locally trivial
elliptic bundle, so that only models $(a)$ and $(d)$ appear. A similar reduction
can be done in the Riccati case, so that only models $(a)$, $(c)$ and $(d)$
appear.

\vskip 2truecm\par
{\bf 8. Poincar\'e metric}
\vglue 1truecm\par

In this section we introduce some analytic tools which will be used to complete
the classification of nongeneral type foliations.
\par
Let ${\cal F}$ be a reduced nef foliation on $X$. Recall that each leaf of
${\cal F}$ is an orbifold, injectively immersed in $X\setminus Sing({\cal F})$.
The universal covering of such a leaf has no multiple
point (indeed, this is true already for the holonomy covering), hence it is
isomorphic to ${\bf D}$ or ${\bf C}$ or ${\bf C}P^1$. In fact, ${\bf C}P^1$ is
not allowed, for the leaf would be a (quotient of a) rational curve $C$ disjoint
from $Sing({\cal F})$ and we would have $K_{\cal F}\cdot C = -\chi_{orb}(C)<0$.
On each leaf of ${\cal F}$ we may therefore put its {\bf Poincar\'e metric}: the
unique complete metric of curvature -1 if the leaf is uniformised by ${\bf D}$
(hyperbolic leaf), and the identically zero ``metric'' if the leaf is 
uniformised by ${\bf C}$ (parabolic leaf).
\par
The following result has been proved in [Br2], and it can be thougth as a
metricised counterpart to Theorem 1. We refer to [Dem] for the basic theory of
singular hermitian metrics.
\par
\underbar{\bf Theorem 5} [Br2]. {\it Suppose that at least one leaf of ${\cal
F}$ is hyperbolic. Then the Poincar\'e metric on the leaves of ${\cal F}$ 
induces on the canonical bundle $K_{\cal F}$ a singular hermitian metric whose 
curvature is a closed positive current.}
\par
More explicitely, this means the following. Let $U$ be an open subset of
$X\setminus Sing(X)$, where ${\cal F}$ is generated by a holomorphic vector
field $v$ with isolated zeroes. This vector field induces a local trivialisation
of $K_{\cal F}$ over $U$. On $U\setminus Sing({\cal F})$ set
$$F=\log \Vert v\Vert_{Poin}$$
where, for each $p\in U\setminus Sing({\cal F})$, $\Vert v(p)\Vert_{Poin}$ is
the Poincar\'e norm of $v(p)$ with respect to the Poincar\'e metric on the leaf
$L_p$ through $p$. Thus $F$ is a function with values into $[-\infty ,+\infty
)$, and $F(p)=-\infty$ iff $L_p$ is parabolic. Moreover, it is
not difficult to see that $F$ is upper semicontinuous. Now, Theorem 5 says that
$F$ is a {\it plurisubharmonic} function. Such a function $F$ is the local
weight of a singular hermitian metric on $K_{\cal F}$, in the local
trivialisation induced by $v$. The curvature $\Omega$ of this metric is locally
expressed by 
$$\Omega ={i\over 2\pi}\partial\bar\partial F$$
and the plurisubharmonicity of $F$ is equivalent to the fact that $\Omega$ is a
{\it closed positive current}. All of this can be done, of course, even if $U$
cuts $Sing(X)$: replace $K_{\cal F}$ by a locally free power of it, etc.
Finally, the fact that $F$ above is plurisubharmonic permits to extend it, in a
plurisubharmonic way, to $Sing({\cal F})$; it turns out that $F=-\infty$ on
$Sing({\cal F})$. Thus the above singular metric on $K_{\cal F}$ is everywhere
defined. The real Chern class of $K_{\cal F}$ is then equal to the De Rham
cohomology class of $\Omega$. 
\par
A first consequence of the theorem is that the existence of a hyperbolic leaf
implies that ``most'' leaves are hyperbolic: parabolic leaves constitute a
pluripolar subset of $X\setminus Sing({\cal F})$. Of course, when all the leaves
are parabolic (we shall say that ${\cal F}$ is a {\bf parabolic foliation}) then
the Poincar\'e metric is not very helpful. We shall see later how to handle with
these parabolic foliations [MQ1] [MQ2].
\par
Another consequence is that the existence of a hyperbolic leaf implies that $\nu
({\cal F})\ge 1$: the current $\Omega$ is strictly positive along hyperbolic
leaves, and so $K_{\cal F}$ is certainly not numerically trivial. 
\par
Let us now discuss the principle of the proof of Theorem 5. Take, in $X\setminus
\{ Sing({\cal F})\cup Sing(X)\}$, a disc $T$ transverse to the foliation. For
each $t\in T$, let $L_t$ be the leaf through $t$, with universal covering
$\widetilde L_t\simeq {\bf D}$ or ${\bf C}$. 
Thanks to the nefness of $K_{\cal F}$
(and also the projectivity, or better the K\"ahlerianity, of $X$) it turns out
that these universal coverings, for $t\in T$, glue together to give a smooth
complex surface $U_T$ [Br2, Proposition 1], called {\bf covering tube}. See also
[Ily]. Such a surface is equipped with a fibration over $T$, $U_T\buildrel
P\over\rightarrow T$, a section $T\buildrel s\over\rightarrow U_T$, 
and an immersion $U_T\buildrel\pi\over\rightarrow X$. For
each $t\in T$ the pointed fibre $(P^{-1}(t),s(t))$ is isomorphic to 
$(\widetilde L_t,t)$, the universal covering of $L_t$ with basepoint $t$, and
$\pi$ sends $P^{-1}(t)$ onto $L_t$ and $s(t)$ to $t$, as a regular (universal)
covering. By construction, the plurisubharmonicity of the
leafwise Poincar\'e metric on $X$ is equivalent to the plurisubharmonicity of
the fibrewise Poincar\'e metric on $U_T$, for every $T$. By results of
Yamaguchi [Yam], this last one follows from a sort of ``holomorphic convexity''
of $U_T$. Thus, the core of [Br2] consists exactly in proving that, for each
$T$, the covering tube $U_T$ possesses this sort of ``holomorphic convexity''.
\par
Another very important property of the leafwise Poincar\'e metric is its
continuity (in the sense that the functions $e^F$ are continuous).
This is one of the major results of [MQ1], and it is based on some
ideas previously introduced in [MQ2].
\par
\underbar{\bf Theorem 6} [MQ1] [MQ2]. {\it The Poincar\'e metric on the leaves
of ${\cal F}$ is continuous. Moreover, the polar set of this metric (i.e.
$Sing({\cal F})\cup\{$parabolic leaves$\}$) is either the full $X$ or a proper
algebraic subset of $X$.}
\par
As for Theorem 5, we discuss here only some of the ideas appearing in the 
proof of Theorem 6.
\par
We may assume that ${\cal F}$ is not a parabolic foliation (otherwise the result
is trivial) nor a fibration (otherwise the result is easy). By a theorem of
Jouanolou [Br1, page 84], ${\cal F}$ has a finite number of algebraic leaves.
Take $p\in X\setminus Sing({\cal F})$ such that the leaf $L_p$ is
transcendental. We shall prove that $L_p$ is hyperbolic and that the Poincar\'e
metric is continuous at $p$.
\par
Suppose, by contradiction, that $L_p$ is parabolic or that the metric is
discontinuous at $p$. In both cases, we may find a sequence $p_n\to p$ and a
sequence of holomorphic maps ${\bf D}\buildrel f_n\over\rightarrow X$ with
$f_n(0)=p_n$ and $f_n({\bf D})\subset L_{p_n}$ such that $\{ f_n\}$ does {\it
not} converge (up to subsequences) to
a map ${\bf D}\buildrel f\over\rightarrow X$ with $f(0)=p$ and
$f({\bf D})\subset L_p$ (see [Br2, Proposition 2]; if $L_p$ is parabolic we may
just choose $p_n=p$ for every $n$ and $f_n(z)=g(nz)$ for some nonconstant 
entire $g:{\bf C}\to L_p$, $g(0)=p$). More precisely [Br2, Proposition 3], there
exists some $r\in (0,1)$ such that $area(f_n({\bf D}_r))\to +\infty$, where
${\bf D}_r\subset {\bf D}$ is the disc of radius $r$ and $area(\ \cdot\ )$ is
computed with respect to any hermitian metric $\omega$ on $X$. This allows to
construct a nontrivial closed positive current $\Phi\in A^{1,1}(X)'$, 
by a variation on Ahlfors' lemma [MQ1, Lemma V.2.5]: 
if $\eta\in A^2(X)$, we set
$$\Phi (\eta )=\lim_{n_k\to +\infty} {\int_0^R{dt\over t}\int_{{\bf D}_t}
f_{n_k}^*(\eta ) \over \int_0^R{dt\over t}\int_{{\bf D}_t}
f_{n_k}^*(\omega )}$$
and it turns out that for some $R\in (r,1)$ and some subsequence $n_k\to
+\infty$ this indeed defines a {\it closed} positive current.
\par
The current $\Phi$ is ${\cal F}$-invariant, and using the fact that $L_p$ is
transcendental (which implies, among other things, that the cohomology class 
$[\Phi ]$ of $\Phi$ is nef) one can prove the following two inequalities:
\par i) $K_{\cal F}\cdot [\Phi ]\le 0$ [MQ1] [MQ2]
\par ii) $N_{\cal F}^*\cdot [\Phi ]\le 0$ [Br5].
\par\noindent 
By Hodge index theorem, the first inequality implies $\nu({\cal F})\le 1$. In
fact we must have $\nu ({\cal F})=1$, otherwise ${\cal F}$ would be parabolic by
Theorem 2 (or Theorem 5). Still by Hodge theorem, we therefore have $c_1(K_{\cal
F})=\lambda [\Phi ]$ for some $\lambda >0$, hence the second inequality becomes
$N_{\cal F}^*\cdot K_{\cal F}\le 0$, i.e. $K_X\cdot K_{\cal F}\le 0$. Now
Riemann--Roch formula gives $h^0(K_{\cal F}^{\otimes n})\ge 1-h^1({\cal O}_X)$
for $n$ large.
Hence if $h^1({\cal O}_X)=0$ then $kod({\cal F})\ge 0$, but the same holds even
if $h^1({\cal O}_X)>0$, by the usual play with the Albanese map (as in Theorems
2 and 3). By Theorem 3 we therefore have $kod({\cal F})=1$. The classification
given by Theorem 4 says that ${\cal F}$ is Riccati or Turbulent. However, in
both cases the Poincar\'e metric on the leaves arises from the pull-back of a
metric on the base $B$ under the rational or elliptic fibration
$X\buildrel\pi\over\rightarrow B$, and in both cases the validity of Theorem 6
is easily checked.
\par
Remarks that the same reasoning gives the classification of parabolic
foliations: they have Kodaira dimension 0 or 1.
\par
We conclude this section by noting that this way to prove the continuity of the
leafwise Poincar\'e metric is not totally satisfying, because it is based on
previous classification results (Theorems 2, 3 and 4). 
It would be better to find an {\it a priori} argument. 
When there are no parabolic leaves at all then this
is indeed possible, using Brody's lemma [Br2, Proposition 2] in place of the
above current $\Phi$.

\vskip 2truecm\par
{\bf 9. Hilbert modular foliations}
\vglue 1truecm\par

With a negligible abuse of terminology, we shall say that a complex projective
surface $X$ (with cyclic quotient singularities) is a {\bf Hilbert modular
surface} [BPV, page 177] if there exists a (possibly empty) curve $C\subset
X\setminus Sing(X)$ such that:
\par\noindent
i) each connected component of $C$ is a cycle of smooth rational curves, 
contractible to a normal singularity (a ``cusp'');
\par\noindent
ii) $X\setminus C$ is uniformised (in orbifold's sense) by the bidisc ${\bf
D}\times{\bf D}$, more precisely $X={\bf D}\times{\bf D}/\Gamma$ where $\Gamma$
is a discrete subgroup of $Aut_0({\bf D}\times{\bf D})=$ $PSL(2,{\bf R})\times
PSL(2,{\bf R})$;
\par\noindent
iii) $\Gamma$ is irreducible, in the sense that it does not contain a finite
index subgroup of the form $\Gamma_1\times\Gamma_2$, with $\Gamma_j\subset
PSL(2,{\bf R})$, $j=1,2$.
\par
Such a surface is naturally equipped with two natural foliations ${\cal F}$ and
${\cal G}$, arising from the horizontal and the vertical foliations on 
${\bf D}\times{\bf D}$, preserved by $\Gamma$. These are called {\bf Hilbert
modular foliations}. They are both tangent to $C$, and they are singular only in
correspondence of the normal crossings of $C$; these singularities are reduced
and nondegenerate (eigenvalues can by computed by Camacho--Sad formula).
\par
We have $K_{\cal F}=N_{\cal G}^*\otimes{\cal O}_X(C)$, the logarithmic conormal
bundle of ${\cal G}$: indeed, $C$ is the tangency locus between ${\cal F}$ and
${\cal G}$, and along $C$ the two foliations have a first order tangency. The
results of section 4 show that $K_{\cal F}$ is nef. By the irreducibility
hypothesis on $\Gamma$ the foliations are not fibrations, therefore by the
logarithmic Castelnuovo--De Franchis--Bogomolov lemma [Br1, Chapter 6] we have
$kod({\cal F})=kod(N_{\cal G}^*\otimes{\cal O}_X(C))\le 0$. The case $kod({\cal
F})=0$ is however excluded because ${\cal F}$ is not a parabolic foliation, so
we finally obtain: ${\cal F}$ is a reduced nef foliation with $\nu ({\cal F})=1$
and $kod({\cal F})=-\infty$. The same of course holds for ${\cal G}$.
\par
The following result completes the classification of nongeneral type foliations.
\par
\underbar{\bf Theorem 7} [Br2] [MQ1]. {\it Let ${\cal F}$ be a reduced nef
foliation with $$\nu ({\cal F})=1\qquad kod({\cal F})=-\infty$$
on a projective surface $X$. Then ${\cal F}$ is a Hilbert modular foliation.}
\par
{\it Proof.}
\par
We give just a sketch, refering to [Br2] for more details.  
\par
The arguments explained around the proof of Theorem 6 show that ${\cal F}$ is
not a parabolic foliation, and so the leafwise Poincar\'e metric is nontrivial.
By Theorem 5, the curvature $\Omega$ of this metric is a closed positive
current. Using $kod({\cal F})=-\infty$ one finds that $\Omega$ has zero Lelong
number everywhere, and even more [Br2, Proposition 5]: $\Omega$ is absolutely
continuous, i.e. its coefficients (measures) are $L^1_{loc}$ functions. This
permits to define the wedge product $\Omega\wedge\Omega$ (in a punctual way, as
a $(2,2)$-form with $L_{loc}^{{1\over 2}}$ coefficients), and from $c_1(K_{\cal
F})=[\Omega ]$, $c_1^2(K_{\cal F})=0$, it follows [Dem] that
$$\Omega\wedge\Omega\equiv 0.$$
Roughly speaking, this suggests that $\Omega$ has a nontrivial Kernel, and we
want to prove that this Kernel integrates to a {\it holomorphic} foliation
${\cal G}$, which will be the ``second'' Hilbert modular foliation on $X$.
\par
Take a local chart ${\bf D}\times{\bf D}\subset X\setminus \{ Sing({\cal F})\cup
Sing(X)\}$, disjoint from the parabolic leaves (which are algebraic by Theorem
6), and in which ${\cal F}$ is expressed by the vector field
${\partial\over\partial w}$. In this local trivialisation we have 
$$\Omega ={i\over 2\pi}\partial\bar\partial F$$ 
where $F$ is finite, continuous and
plurisubharmonic, by Theorems 5 and 6. Moreover, the second derivatives
$F_{\alpha\bar\beta }$ ($\alpha ,\beta\in\{ z,w\}$) are 
$L_{loc}^1$ ($\Omega$ absolutely
continuous), $F_{z\bar z}F_{w\bar w}-F_{z\bar w}F_{w\bar z}$ is identically zero
($\Omega\wedge\Omega\equiv 0$), and $F_{w\bar w}=e^F$ (curvature -1 along the
leaves). The Kernel of $\Omega$ is expressed by the differential equation
$${dw\over dz}=-{F_{z\bar w}\over F_{w\bar w}} $$
and we shall see below (Proposition 1, [Br2, Proposition 6]) that the function 
$-F_{z\bar w}/F_{w\bar w}$ is {\it holomorphic}. 
By integrating the differential equation,
we obtain, at least outside the parabolic leaves and the singularities, a
holomorphic foliation ${\cal G}$, transverse to ${\cal F}$. Using again the 
absolute continuity of $\Omega$, one proves that ${\cal G}$ extends to the full
$X$, as a holomorphic foliation with a first order tangency with ${\cal F}$ 
along the parabolic leaves: the above function $-F_{z\bar w}/F_{w\bar w}$
extends to a meromorphic function with a first order pole along $\{
F=-\infty\}$. We thus have $N_{\cal G}^*\otimes{\cal O}_X(C)=K_{\cal F}$, where
$C$ is the closure of the parabolic leaves.
\par
The next step consists in analysing the structure of $C$. By numerical arguments
similar to those employed in the step 2) in the proof of Theorem 3, we obtain
that each connected component of $C$ is a contractible cycle of rational 
curves, and ${\cal F}$ and ${\cal G}$ 
are singular only in correspondence of its normal crossings.
The singularities of ${\cal G}$ are nondegenerate, and using Camacho--Sad
formula we see also that they are reduced. Finally, by arguments close to those
before the statement of Theorem 7 ($K_{\cal G}=N_{\cal F}^*\otimes{\cal
O}_X(C)$) we obtain that ${\cal G}$ is nef and its Kodaira dimension is
$-\infty$.
\par 
Therefore we may repeat the previous steps, for
${\cal G}$ instead of ${\cal F}$, and we find a third holomorphic foliation
${\cal H}$, generated by the Kernel of the curvature of the Poincar\'e metric on
the leaves of ${\cal G}$. Without surprise, one proves that ${\cal H}$ 
coincides with ${\cal F}$. 
\par
Thus, outside the parabolic leaves of ${\cal F}$ (which
coincide with those of ${\cal G}$, and which form cycles $C$ of rational 
curves) we have a pair of transverse foliations such
that the leaves of one are in the Kernel of the curvature of the Poincar\'e
metric on the leaves of the other. This means that the holonomy maps between
pieces of leaves of ${\cal F}$ (resp. ${\cal G}$) induced by ${\cal G}$ (resp.
${\cal F}$) are {\it isometries} in the Poincar\'e metric (because they are
holomorphic and they preserve the Poincar\'e area form). As a consequence of
this, we can define a complete hermitian metric on $X\setminus C$ which is
locally isometric to the Poincar\'e metric on ${\bf D}\times{\bf D}$, by taking
the orthogonal sum of the metrics along ${\cal F}$ and along ${\cal G}$. It is
then easy to conclude, by classical arguments, that $X\setminus C$ is 
uniformised by ${\bf D}\times{\bf D}$, and then that ${\cal F}$ 
and ${\cal G}$ are Hilbert modular foliations.
\par
$\triangle$
\par
In the course of the proof we have used the following quite miraculous fact
concerning Monge--Amp\`ere foliations.
\par
\underbar{\bf Proposition 1}. {\it Let $F:{\bf D}\times{\bf D}\rightarrow {\bf
R}$ be a continuous plurisubharmonic function such that:
\par\noindent
i) $\Omega ={i\over 2\pi}\partial\bar\partial F$ is absolutely continuous;
\par\noindent
ii) $\Omega\wedge\Omega\equiv 0$;
\par\noindent
iii) $F_{w\bar w}=e^F$.
\par\noindent
Then $-F_{z\bar w}/F_{w\bar w}$ (the slope of the Kernel of $\Omega$) is a
holomorphic function.}
\par
{\it Proof}.
\par
A full proof is given in [Br2, Proposition 6], here we shall explain a
different, more geometric, proof, which however requires some additional
regularity. More precisely, we shall suppose here that $F$ is {\it smooth} 
(using standard elliptic regularity theory one can prove that, 
under the hypotheses of the proposition, 
$F$ Lipschitz implies $F$ smooth; however, the step from
continuous to Lipschitz is more problematic).
\par
The Kernel of $\Omega$ being smooth and $\Omega$ being closed, we certainly have
a smooth foliation ${\cal G}$ with complex leaves obtained by integrating
that Kernel, a so-called Monge--Amp\`ere foliation [Kli]. Generally speaking, 
such a foliation is far from being holomorphic, and we need to prove 
that indeed it is, thanks to hypothesis iii).
\par
Let $L_0$ be a leaf of ${\cal G}$, say the one through $(0,0)$; we may suppose
$L_0=\{ w=0\}$. The leaf $L_t$ through $(0,t)$, $t$ close to 0, has equation
$$w=a(z)t+b(z)\bar t + o(\vert t\vert )$$
where $a,b$ are holomorphic functions and $a(0)=1$, $b(0)=0$. Let 
$$\phi_s : \{ z=0\} \rightarrow \{ z=s\}$$
be the holonomy map defined by ${\cal G}$, on a neighbourhood of $(0,0)\in \{
z=0\}$. We have
$$\phi_s^*(dw\wedge d\bar w)=(\vert a(s)\vert^2-\vert b(s)\vert^2)dw\wedge d\bar
w + O(\vert w\vert )$$
and so, by hypothesis iii),
$$\phi_s^*(\Omega\vert_{\{ z=s\} }) = 
e^{F(s,0)-F(0,0)}(\vert a(s)\vert^2-\vert b(s)\vert^2)\Omega\vert_{\{ z=0\} } 
+ O(\vert w\vert ) .$$
But ${\cal G}$ is in the Kernel of $\Omega$ and $\Omega$ is closed, hence
$\phi_s^*(\Omega\vert_{\{ z=s\} })$ must be equal to $\Omega\vert_{\{ z=0\} }$
by Stokes theorem. In particular:
$$\log (\vert a(s)\vert^2-\vert b(s)\vert^2) = F(s,0)-F(0,0).$$
The function $F(s,0)$ is harmonic, and not only subharmonic, because
$\Omega\vert_{L_0}\equiv 0$. Thus $\log (\vert a(s)\vert^2-\vert b(s)\vert^2)$
is harmonic too, and the maximum principle plus the initial conditions $a(0)=1$,
$b(0)=0$ give $b\equiv 0$. In other words, the linear part at 0 of the holonomy
diffeomorphism $\phi_s$ is a complex linear map, for every $s$, and so the
foliation ${\cal G}$ is holomorphic along $L_0$. This leaf being an arbitrary
leaf, we obtain that ${\cal G}$ is holomorphic {\it tout court}.
\par
$\triangle$

\vskip 2truecm\par
{\bf 10. K\"ahler surfaces}
\vglue 1truecm\par

Up to now we have mostly considered the case of foliations on projective
surfaces. Let us explain in this last section how to extend the previous results
to the slightly more general case of compact K\"ahler surfaces. More precisely,
we shall draw a complete list of {\it all} foliations on compact K\"ahler 
nonprojective surfaces: this is consistent with the fact that on such a surface
a foliation cannot be of general type, because a line bundle cannot have Kodaira
dimension 2.
\par
Instead of working with a special model of the foliation (reduced, nef), it  
is simpler here to work with a special model of the surface. Thus, let $X$ be a
{\it minimal} smooth compact K\"ahler nonprojective surface. 
According to Kodaira classification [BPV], the algebraic dimension 
$a(X)$ is either 1, in which case $X$ has a unique elliptic fibration, or 0, 
in which case $X$ is a torus or a K3 surface. 
Let ${\cal F}$ be any foliation on $X$. We distinguish three cases:
\par
\underbar{$a(X)=1$}. Let $\pi :X\to B$ be the elliptic fibration. Any compact
curve in $X$ is contained in some fibre of $\pi$. In particular, if ${\cal F}$
is different from the fibration then the tangency locus between ${\cal F}$ and
$\pi$ is contained in some fibres. In other words, ${\cal F}$ is transverse to
a generic fibre of $\pi$, and so ${\cal F}$ is a {\bf Turbulent foliation}.
\par
\underbar{$X$ is a torus with $a(X)=0$}. We may take on $X={\bf C}^2/G$ a
foliation ${\cal G}$ generated by a constant vector field and tangent to ${\cal
F}$ at some point. The tangency locus between ${\cal F}$ and ${\cal G}$ is not
empty and it cannot be a curve, because such an $X$ contains no curve. Thus it
remains only the possibility ${\cal F}={\cal G}$, i.e. ${\cal F}$ is a {\bf
Kronecker foliation}, generated by a constant vector field.
\par
\underbar{$X$ is a K3 surface with $a(X)=0$}. Such a surface admits, by Yau's
theorem [BPV, page 40], 
a K\"ahler--Einstein metric in each K\"ahler cohomology class, which
implies the semistability of its tangent bundle [Fri, Chapter 4]. 
That is, if $\omega\in
A^{1,1}(X)$ is any K\"ahler form on $X$, then we have the semistability
inequality
$$c_1(T_{\cal F})\cdot [\omega ]\le {1\over 2}c_1(X)\cdot [\omega ] =0$$
i.e.
$$c_1(K_{\cal F})\cdot [\omega ]\ge 0.$$
This inequality plays here the same r\^ole as Theorem 1 in the projective case:
it says that $K_{\cal F}$ is pseudoeffective [Dem] [Lam]. Looking at the proof
of [Fuj], one realizes that such a line bundle has a Zariski decomposition: as a
${\bf Q}$-bundle, it can be written as $$K_{\cal F}=P+N$$ where $N$ is an
effective ${\bf Q}$-bundle with contractible support and $P$ is a ${\bf
Q}$-bundle which satisfies: i) $P\cdot [\omega ]\ge 0$ for every K\"ahler form
$\omega$; ii) $P\cdot C\ge 0$ for every compact curve $C\subset X$; iii) $P\cdot
D=0$ for every compact curve $D\subset Supp(N)$. From i) and ii) it follows 
[Lam] that $P\cdot P\ge 0$. 
But on a K3 surface with zero algebraic dimension there exists
only one line bundle $L$ with $L\cdot L\ge 0$: the trivial one (by
Riemann--Roch we have $h^0(L)+h^0(L^*)\ge 2$ and by $a(X)=0$ we have $h^0(L)\le
1$ and $h^0(L^*)\le 1$). It follows that $P$ is numerically trivial, or more
precisely $K_{\cal F}^{\otimes n}={\cal O}_X(E)$ for some $n>0$ and some
effective divisor $E$ supported on $Supp(N)$.
\par
We can now repeat the arguments of sections 4 and 5 (with some care, because
${\cal F}$ could be nonreduced). After contraction of $Supp(N)$ we obtain a
foliation ${\cal F}'$ on a surface $X'$ with cyclic quotient singularities,
whose canonical bundle $K_{{\cal F}'}$ is a torsion ${\bf Q}$-bundle. 
After a regular
covering $Y\buildrel\pi\over\rightarrow X'$ this canonical bundle becomes
holomorphically trivial. Of course $Y$ has still algebraic dimension 0, and
$Y\not= X'$, because a K3
surface has no global holomorphic vector field. It follows that $Y$ is a torus
and $\pi^*({\cal F}')$ is a Kronecker foliation. Thus ${\cal
F}'$ is a {\bf Kummer foliation}, quotient of a Kronecker foliation by a
$n$-cyclic group.
 
\vskip 2truecm\par
\centerline{\bf References}
\vglue 1truecm\par\noindent
[And] A. Andreotti, {\it Th\'eor\`emes de d\'ependance alg\'ebrique sur les
espaces complexes pseudoconcaves}, Bull. Soc. Math. France 91 (1963), 1-38
\par\noindent
[Ara] S. Arakelov, {\it Families of algebraic curves with fixed degeneracy},
Math. USSR Izv. 5 (1971), 1277-1302
\par\noindent 
[B-M] F. Bogomolov, M. McQuillan, {\it Rational curves on foliated varieties},
preprint IHES M/01/07 (2001)
\par\noindent
[Br1] M. Brunella, {\it Birational geometry of foliations}, First Latin American
Congress of Mathematicians, Notas de Curso, IMPA (2000) (available also at
www.impa.br)
\par\noindent
[Br2] M. Brunella, {\it Subharmonic variation of the leafwise Poincar\'e
metric}, Inv. Math. (to appear)
\par\noindent
[Br3] M. Brunella, {\it Minimal models of foliated algebraic surfaces}, Bull.
Soc. Math. France 127 (1999), 289-305
\par\noindent
[Br4] M. Brunella, {\it Complete polynomial vector fields on the complex plane},
preprint Univ. Bourgogne 290 (2002)
\par\noindent
[Br5] M. Brunella, {\it Courbes enti\`eres et feuilletages holomorphes}, L'Ens.
Math. 45 (1999), 195-216
\par\noindent
[C-S] C. Camacho, P. Sad, {\it Invariant varieties through singularities of
holomorphic vector fields}, Ann. Math. 115 (1982), 579-595
\par\noindent
[C-M] D. Cerveau, R. Moussu, {\it Groupes d'automorphismes de $({\bf C},0)$ et
\'equations diff\'erentielles $ydy+...=0$}, Bull. Soc. Math. France 116 (1988),
459-488
\par\noindent
[Dem] J.-P. Demailly, {\it $L^2$-vanishing theorems for positive line bundles
and adjunction theory}, Transcendental Methods in Algebraic Geometry, Springer
Lecture Notes 1646 (1996), 1-97
\par\noindent
[Fri] R. Friedman, {\it Algebraic surfaces and holomorphic vector bundles},
Springer (1998)
\par\noindent
[Fuj] T. Fujita, {\it On Zariski problem}, Proc. Japan Acad. A 55 (1979),
106-110
\par\noindent
[Ily] Ju. Il'yashenko, {\it Covering manifolds for analytic families of leaves
of foliations by analytic curves}, Topol. Meth. Nonlin. Anal. 11 (1998), 361-373
\par\noindent
[Kli] M. Klimek, {\it Pluripotential theory}, Clarendon Press, Oxford (1991)
\par\noindent
[Lam] A. Lamari, {\it Le c\^one K\"ahl\'erien d'une surface}, J. Math. Pures
Appl. 78 (1999), 249-263
\par\noindent
[L-M] F. Loray, R. Meziani, {\it Classification de certains feuilletages
associ\'es \`a un cusp}, Bol. Soc. Bras. Mat. 25 (1994), 93-106
\par\noindent
[MQ1] M. McQuillan, {\it Noncommutative Mori theory}, preprint IHES M/00/15
(2000) (revised: M/01/42 (2001))
\par\noindent
[MQ2] M. McQuillan, {\it Diophantine approximations and foliations}, Publ. 
Math. IHES 87 (1998), 121-174
\par\noindent
[MR1] J. Martinet, J.-P. Ramis, {\it Probl\`emes de modules pour des \'equations
diff\'erentielles non lin\'eaires du premier ordre}, Publ. Math. IHES 55
(1982), 63-124
\par\noindent
[MR2] J. Martinet, J.-P. Ramis, {\it Classification analytique des \'equations
diff\'erentielles non lin\'eaires r\'esonnantes du premier ordre}, Ann. Sci. ENS
16 (1983), 571-621
\par\noindent 
[M-M] J.-F. Mattei, R. Moussu, {\it Holonomie et int\'egrales premi\`eres}, Ann.
Sci. ENS 13 (1980), 469-523
\par\noindent
[Men] L. G. Mendes, {\it Kodaira dimension of singular holomorphic foliations},
Bol. Soc. Bras. Mat. 31 (2000), 127-143
\par\noindent
[Miy] Y. Miyaoka, {\it Deformation of a morphism along a foliation and
applications}, Algebraic Geometry, Bowdoin 1985, Proc. Symp. Pure Math. 46
(1987), 245-268
\par\noindent
[M-P] Y. Miyaoka, Th. Peternell, {\it Geometry of higher dimensional algebraic
varieties}, Birkh\"auser (1997)
\par\noindent
[Pau] E. Paul, {\it Feuilletages holomorphes \`a holonomie r\'esoluble}, J.
reine angew. Math. 514 (1999), 9-70
\par\noindent
[Sak] F. Sakai, {\it Weil divisors on normal surfaces}, Duke Math. J. 51 (1984),
877-887
\par\noindent
[Ser] F. Serrano, {\it Fibred surfaces and moduli}, Duke Math. J. 67 (1992),
407-421
\par\noindent
[ShB] N. I. Shepherd-Barron, {\it Miyaoka's theorems on the generic
seminegativity of $T_X$ and on the Kodaira dimension of minimal regular
threefolds}, Flips and abundance for algebraic threefolds, Ast\'erisque 211
(1992), 103-114
\par\noindent
[Suw] T. Suwa, {\it Indices of vector fields and residues of holomorphic
singular foliations}, Hermann (1998)
\par\noindent
[Szp] L. Szpiro, {\it Propri\'et\'es num\'eriques du faisceau dualisant
relatif}, Pinceaux de courbes de genre au moins deux, Ast\'erisque 86 (1981),
44-78
\par\noindent
[Yam] H. Yamaguchi, {\it Calcul des variations analytiques}, Jap. J. Math. 7
(1981), 319-377

\end